\documentclass{article}
\usepackage{subcaption}
\usepackage[utf8]{inputenc}
\usepackage{natbib}
\usepackage[table,xcdraw]{xcolor}
\usepackage{booktabs}
\usepackage{graphicx}
\usepackage{amsmath}
\usepackage{amsfonts}
\usepackage{color}
\usepackage{algorithm}
\usepackage{algpseudocode}
\usepackage{comment}
\algnewcommand\algorithmicswitch{\textbf{switch}}
\algnewcommand\algorithmiccase{\textbf{case}}
\algnewcommand\algorithmicassert{\texttt{assert}}
\algnewcommand\Assert[1]{\State \algorithmicassert(#1)}%
\algdef{SE}[SWITCH]{Switch}{EndSwitch}[1]{\algorithmicswitch\ #1\ \algorithmicdo}{\algorithmicend\ \algorithmicswitch}%
\algdef{SE}[CASE]{Case}{EndCase}[1]{\algorithmiccase\ #1}{\algorithmicend\ \algorithmiccase}%
\algtext*{EndSwitch}%
\algtext*{EndCase}%

\title{A New Class of Compact Formulations for Vehicle Routing Problems}
\author{Udayan Mandal\textsuperscript{ \rm 1}, Amelia Regan\textsuperscript{ \rm 2}, Louis Martin Rousseau \textsuperscript{\rm 3}, Julian Yarkony \textsuperscript{\rm 4} \\
\textsuperscript{\rm 1} Stanford University, Palo Alto, California USA\\
\textsuperscript{\rm 2} University of Washington, Seattle, Washington USA\\
\textsuperscript{\rm 3} Polytechnique Montreal, Montreal, Quebec, Canada\\
\textsuperscript{\rm 4} Optym
}
\date{February, 2024}

\begin{document}

\maketitle

\begin{abstract}
This paper introduces a novel compact mixed integer linear programming (MILP) formulation and a discretization discovery-based solution approach for the Vehicle Routing Problem with Time Windows (VRPTW). We aim to solve the optimization problem efficiently by constraining the linear programming (LP) solutions to use only flows corresponding to time and capacity-feasible routes that are locally elementary (prohibiting cycles of customers localized in space).

We employ a discretization discovery algorithm to refine the LP relaxation iteratively. This iterative process alternates between two steps: (1) increasing time/capacity/elementarity enforcement to increase the LP objective, albeit at the expense of increased complexity (more variables and constraints), and (2) decreasing enforcement without decreasing the LP objective to reduce complexity. This iterative approach ensures we produce an LP relaxation that closely approximates the optimal MILP objective with minimal complexity, facilitating an efficient solution via an off-the-shelf MILP solver.

The effectiveness of our method is demonstrated through empirical evaluations on classical VRPTW instances. We showcase the efficiency of solving the final MILP and multiple iterations of LP relaxations, highlighting the decreased integrality gap of the final LP relaxation. We believe that our approach holds promise for addressing a wide range of routing problems within and beyond the VRPTW domain.
\end{abstract}
\section{Introduction}
This paper proposes a novel compact formulation for the Vehicle Routing Problem with Time Windows (VRPTW) with an empirically ``tighter" linear programming (LP) relaxation compared to a standard two-index compact formulation baseline, which can be solved directly as a Mixed Integer Linear Program (MILP). By ``tighter," we mean that the LP relaxation's optimal solution objective value is closer to the MILP objective value relative to a baseline LP. We characterize such a formulation with two properties: sufficiency and parsimony. 
Sufficiency indicates that our formulation describes the tightest LP within a considered class of LPs, while parsimony ensures that our formulation is as compact as possible (measured by the number of variables and constraints) without decreasing the LP objective. Utilizing these properties, we devise an algorithm that iteratively expands and contracts this formulation and demonstrates its monotonic convergence to a solution that is both sufficient and parsimonious. By producing a sufficient and parsimonious MILP we aim to ensure that the MILP can be solved efficiently. Moreover, adjusting the formulation parameters allows us to explore the trade-off between solution time and quality.

Our approach augments the standard compact formulation by enforcing agreement between the standard compact decision variables and separate sets of variables and constraints designed to enforce elementarity, time feasibility, and capacity feasibility. Agreement between these formulations and the compact formulation requires that the number of edges indicating travel from one location to another agrees with the original compact variables. Restrictions on time and capacity feasibility are enforced using discretization of time/capacity ranges into ``buckets" , inspired by previous work \citep{boland2017continuous,boland2019price}, while restrictions on elementarity use Local Area Arcs from our recent work \citep{mandal2023graph,mandal2022local_2}. We elaborate on our use of these enforcement mechanisms below. 

\textbf{Restricting Violations of Capacity/Time Feasibility:} We partition the remaining capacity before servicing each customer into ``buckets", creating nodes for a directed acyclic graph. Each pair of nodes corresponding to feasible flow is connected and associated with a decision variable. A valid solution associated with this directed acyclic graph is a weighted non-negative sum of source-to-sink paths.  Operations aiming towards sufficiency reduce the granularity of the partitions. Operations aiming towards parsimony increase the granularity without decreasing the LP objective.  This approach is used for time-enforcing temporal feasibility as well. 

\textbf{Restricting Violations of Elementarity: } For each customer, we create variables representing all possible immediately succeeding sequences of nearby customers. These sequences are limited to those that can be used in an optimal solution, given the intermediate customers, those on the efficient frontier of minimal cost, earliest completion time, and latest start time. We require one variable to be selected for each customer. Operations aiming toward sufficiency expand the set of nearby customers and, hence, the number of variables. This limits the set of cycles of customers allowed, at the cost of increasing the complexity of the LP. Operations aiming towards parsimony analyze the dual solution to decrease the set of nearby customers, decreasing the complexity of the LP in such a manner as not to decrease the LP objective.  We should note that the use LA-arcs for cycle elimination, enables users to expand their formulation's scale as MILP solvers' computational capabilities advance.

We organize this document as follows. Section \ref{sec_litRev} briefly reviews the most relevant literature. Section \ref{sec_tech_review} presents the VRPTW and its compact two-index formulation. Section \ref{sec_much_tighter} provides the proposed LP relaxation for the VRPTW. Section \ref{sec_expand_contract} provides an algorithm called LA-Discretization to solve the corresponding relaxation efficiently. Section \ref{my_r_comp} we consider LA-Arc generation in detail. Section \ref{sec_exper} provides experimental results comparing LA-Discretization to a baseline compact two-index formulation. Section \ref{sec_conc} provides conclusions and discusses extensions.  
\section{Literature Review}
\label{sec_litRev}
Different MILP formulations for VRPTW trade off the increase in the number of variables/constraints against the tightness of the LP relaxation.  The use of a tighter formulation often means that branch \& bound expands fewer nodes; however, this tightness often comes at the cost of a larger LP (in terms of more variables/constraints), which is thus more computationally intensive to solve at each node of the branch \& bound solver. Well-known expanded or Dantzig-Wolfe formulations \citep{dantzig1960decomposition} often provide much tighter linear programming (LP) relaxations than a standard compact formulation for MILPs in many application areas \citep{barnprice,geoffrion1974lagrangean}. These formulations are solved through Column Generation (CG) \citep{cuttingstock} and are widely used across many problem domains, including vehicle routing \citep{Desrochers1992}, crew scheduling \citep{desaulniers1997crew}, and more recently applications in artificial intelligence domains, such as computer vision \citep{yarkony2020data} and machine learning \citep{FlexDOIArticle}. 

Expanded formulations used to solve Vehicle Routing Problems (VRP) are weighted set covering formulations where each column (primal variable) corresponds to a feasible route, and each cover constraint (primal constraint) corresponds to a customer, which must be serviced. The expanded formulation circumvents multiple key difficulties of a compact formulation, the most important of which is the following. The expanded formulation enforces that the flows of the corresponding compact formulation only describe feasible routes.  
However, expanded formulations present multiple difficulties including the following.  
\begin{itemize}
 \item When the LP solution includes fractional variables at the termination of CG, the set of columns (variables in the MILP) generated over the course of CG does not necessarily include all columns that are part of the optimal integer solution.  Thus, we cannot use a standard off-the-shelf MILP solver to guarantee an optimal solution. The importance of this difficulty is emphasized in \citep{barnprice}, which observes that columns that describe an optimal solution to the weighted set cover LP may not express an optimal or even a near-optimal solution to the weighted set cover integer linear program (ILP). Branch \& Price \citep{barnprice} or Branch-Cut \& Price \citep{gauvin2014branch,pecin2017improved} can be used to solve these formulations to optimality, but are difficult to implement. 
 \item Generating negative (or positive in the case of a maximization objective) reduced cost columns is often a challenging optimization problem. In the context of VRPs, these boil down to an elementary resource-constrained shortest path problem \citep{irnich2005shortest} (ERCSPP), which is NP-hard. ERCSPP solvers are tricky to implement, which has led researchers to propose numerous sophisticated implementations in the last decades \citep{righini2008new,righini2009decremental,baldacci2011new,boland2006accelerated,righini2006symmetry,lozano2016exact,tilk2017asymmetry,beasley1989algorithm,righini2004dynamic,cabrera2020exact,lozano2013exact,li2019revised}. 
 \item Valid inequalities such as subset-row inequalities \citep{jepsen2008subset} or rounded capacity inequalities \citep{archetti2011column} must be separated by the CG developer, and these may alter the structure of the CG pricing problem (as is the case for \citep{jepsen2008subset}), which can add further complexity to the pricing problem as additional resources must be included in ERCSPP. 
 \item  CG acceleration methods are often required to ensure fast convergence of CG such as dual stabilization:\citep{du1999stabilized,marsten1975boxstep,haghani2020smooth,rousseau2007interior}, heuristic pricing \citep{desrosiers2005primer}, and complimentary columns \citep{ghoniem2009complementary} 
\end{itemize}

 Elementarity, the property that a customer may be present at most once in a route, can be very challenging to enforce in pricing.  The following approaches attempt to avoid fully/explicitly considering elementarity in routes in pricing problems in CG formulations:   neighborhood (NG)-route relaxations \citep{baldacci2011new,baldacci2012recent}, decremental state space relaxations (DSSR) \citep{righini2009decremental}, P-step formulations \citep{dollevoet2020p}, Q-route relaxations \citep{christofides1981exact} and Local Area (LA)-route relaxations \citep{mandal2022local_2,mandal2023graph}. In NG-route relaxations \citep{baldacci2012recent}, the elementarity of the route is relaxed only to prevent cycles of customers localized in space. Specifically, when extending a label during pricing, we hold onto customers from the previous label that are in the NG-neighborhood (nearby customers) of the latest customer at the given label.  We forbid the extension of a label to a customer where the label states that it has already been. DSSR can be understood as expanding the NG-neighborhood as needed \citep{righini2008new,righini2009decremental}, hence, gradually enforcing elementarity where violated. Q-route relaxations relax elementarity to prevent cycles that contain any sub-sequence in which a customer is preceded and succeeded by the same customer.  P-step formulations \citep{dollevoet2020p} employ an arc-based formulation where arcs correspond to elementary sequences of customers. CG is employed to solve optimization since the number of such sequences can be quite large. LA-routes \citep{mandal2022local_2} construct routes as sequences of elementary sub-sequences of customers where each sub-sequence is localized in space.  

The work of \citep{boland2017continuous} with respect to continuous-time service network design problems can be used to enforce feasibility with regard to time in VRPTWs. In that paper, the authors intelligently construct a partition of the times in a time-expanded graph representation of the underlying problem. This partition induces a relaxation of time feasibility. Specifically, a vehicle can depart at any time in an interval, given that it arrives at any time in that interval. The weakness is that this partition permits a vehicle to leave a node before it arrives. No relaxation occurs if the partition is constructed such that arrivals at nodes only occur at the earliest time possible in the optimal solution. To construct such a partition, the authors iterate between solving the MILP and splitting used nodes that cause a violation of temporal feasibility.

The approach that we propose in this paper is able to solve a large portion of the problems in the Solomon instances though not as many as Branch \& Price based solvers \citep{baldacci2012recent}.  

\section{The Vehicle Routing Problem with Time Windows }
\label{sec_tech_review}
This section reviews the core mathematical concepts required to understand our paper.  We organize this section as follows. In Section \ref{subsec_cvrp_no_notation}, we informally review the Vehicle Routing Problem with Time Windows (VRPTW). Section \ref{subsec_math_notanoi} provides the key notation used to discuss VRPTW, while Section \ref{sec_classic_compact} describes a classic compact two-index formulation for VRPTW.
\subsection{Informal Description}
\label{subsec_cvrp_no_notation}
A VRPTW problem instance is described using the following terms: 
\textbf{(a)} a depot located in space; \textbf{(b)} a set of customers located in space, each of which has an integer demand, service time, and time window for when service starts; and \textbf{(c)} an unlimited number of homogeneous vehicles with integer capacity.  

A solution to VRPTW assigns vehicles to routes, where each route satisfies the following properties: \textbf{(1)}  The route starts and ends at the depot. \textbf{(2)} The total demand of customers serviced on the route does not exceed the vehicle's capacity. \textbf{(3)} The route cost is the total distance traveled. \textbf{(4)} The vehicle leaves a customer immediately after servicing that customer.  
\textbf{(5)} Service starts at a given customer within the time window of that customer. \textbf{(6)} A customer is serviced at most once in a route.  

The VRPTW solver selects a set of routes with the goal of minimizing the total distance traveled while ensuring that each customer is serviced exactly once across all selected routes.
\subsection{Mathematical Notation}
\label{subsec_math_notanoi}
We use $N$ to denote the set of customers, which we index by $u$. We use $N^+$ to denote $N$ augmented with the starting/ending depot (which are co-located), and denoted $\alpha, \bar{\alpha}$ respectively.  
Each customer $u\in N$ has a demand denoted $d_u$. We define the demand at the depots as $0$, which we write formally as $d_{\alpha}=d_{\bar{\alpha}}=0$.  We are given an unlimited number of homogeneous vehicles with capacity $d_0$ that start at the depot at time $t_0$ where $t_0$ is the length of the time horizon of optimization. Note that this is without loss of generality as we can easily apply a cost (distance) term to each departure/arrival from the depot to restrict the number of vehicles deployed. 
We use $t^+_u$ and $t^-_u$ to denote the earliest/latest times for the service at customer $u$ to begin.  The service windows for the start/end depot are denoted by $t_{\alpha}^{+}=t_{\alpha}^{-}=t_{\bar{\alpha}}^{+}=t_0$ and $t_{\bar{\alpha}}^{-}=0$.  Servicing a given customer $u\in N$ takes $t_u^*$ units of time, where $t_{\alpha}^*=t_{{\bar{\alpha}}}^*=0$.

For each $u \in N^+,v\in N^+$ we use $c_{uv}=t_{uv}$ to denote the cost and time required to travel from $u$ to $v$. For simplicity, $t_{uv}$ is altered to include the service time at customer $u$, which is $t_u^*$.  Thus, we set $t_{uv}\leftarrow t_{uv}+t_{u}^*$, then set $c_{uv}\leftarrow t_{uv}$, which offsets the cost of a solution by a constant. 

\subsection{ A Two-Index Compact Formulation}
\label{sec_classic_compact}
We now describe the two-index formulation for the VRPTW of interest.  We use the following decision variables. We set binary decision variable $x_{uv}=1$ to indicate that a vehicle departs $u$ and travels immediately to $v$.  $x_{uv}$ is defined for every pair of customers/depots, where $u$ is not the ending depot and $v$ is not $u$ or the starting depot. The term $x_{uv}$ is also not defined for cases where the route that proceeds from beginning to end $[\alpha,u,v,\bar{\alpha}]$ is infeasible due to time or capacity. The set of existing $x_{uv}$ terms are denoted $E^*$.

We use continuous decision variable $\tau_u \in \mathbb{R}_{0+}$ to denote the amount of time remaining when service starts at customer $u$.  We enforce the bounds on the start of the service with the following inequalities:  $t^+_u\geq \tau_u\geq t^-_u \quad \forall u \in N$.

We use continuous decision variable $\delta_u \in \mathbb{R}_{0+}$ to denote the capacity remaining in the vehicle immediately before servicing customer $u$.  We enforce that sufficient capacity remains to service customer $u$ with the following inequalities $d_0\geq \delta_u\geq d_u \quad \forall u \in N$. 

We now write the two-index form for the VRPTW as a compact MILP :
\begin{subequations}
\label{orig_opt}
    \begin{align}
      &  \min_{\substack{x_{uv} \in \{0,1\} \\ t^+_u\geq \tau_u\geq t^-_u\\ d_0\geq  \delta_u \geq d_u}} 
\sum_{\substack{uv \in E^*}}c_{uv}x_{uv} & \label{comp_obj}\\
    &    \sum_{\substack{ uv \in E^*}}x_{uv}= 1 &\quad \forall u \in N  \quad \label{serv_eq}\\
    &   \sum_{vu \in E^*}x_{vu}=1 &\quad \forall u \in N \label{flow_con}\\
     &   \delta_v-d_v \geq \delta_u-(d_0+d_v)(1-x_{vu}) & \quad \forall u \in N, vu \in E^*  \label{demand_con}\\
      &  \tau_v-t_{vu} \geq \tau_u-(t^+_u+t_{vu})(1-x_{vu}) &\quad \forall u \in N ,vu \in E^*\label{time_con}\\
       & \sum_{u \in N}x_{\alpha u}\geq \lceil \frac{\sum_{u \in N}d_u}{d_0}\rceil& \label{minVeh_con}
        \end{align}
        \end{subequations}
In \eqref{comp_obj} we minimize the total distance traveled by all vehicles. In \eqref{serv_eq}, we enforce that each customer is serviced exactly once. In \eqref{flow_con}, we enforce that a vehicle must leave each customer once. In \eqref{demand_con}, we enforce that if $u$ is immediately preceded by $v$ in a route, then the amount of demand remaining immediately before servicing $u$ (denoted $\delta_u$) is no greater than that of $v$ (denoted $\delta_v$) minus $d_v$; and is otherwise unrestricted. Observe that in \eqref{demand_con} that $(d_0+d_v)(1-x_{vu})$ operates as a big-M term, which makes the inequality inactive when $x_{vu}=0$.  In \eqref{time_con}, we enforce that if $u$ is preceded by $v$ in a route, then the amount of time remaining immediately before servicing $u$ (denoted $\tau_u$) is no greater than that of $v$ (denoted $\tau_v$) minus $t_{vu}$; and is otherwise unrestricted. Observe that in \eqref{time_con} that $(t^+_u+t_{vu})(1-x_{vu})$ operates as a big-M term, which makes the inequality inactive when $x_{vu}=0$.

In \eqref{minVeh_con} we enforce that a minimum number of vehicles exit the depot; equal to a lower bound on the number of vehicles required to service all demand. This lower bound is the sum of the demands divided by the vehicle capacity, then rounded up to the nearest integer. This is an optional constraint that is not common in the CG literature. 

It is well know that this formulation can lead to weak LP bounds because of the following issues. (1)  It does not explicitly eliminate sub-tours in the LP solution. (2)  Capacity and time feasibility are not well enforced.
\section{A Tighter Compact LP Relaxation}
\label{sec_much_tighter}

In this section, we introduce our compact formulation for VRPTW, which is, in our experiments, tighter than the two-index compact formulation in \eqref{orig_opt}. This formulation is parameterized by sets constructed later in the document (see Section \ref{sec_expand_contract}). Specifically, our formulation is parameterized by unique values of the following for each customer: \textbf{(1)} set of Local Area (LA) neighbors which limit the class of cycles permitted in LP solutions \citep{mandal2022local_2} (which we discuss in Section \ref{subsec_la_incor}), \textbf{(2)} a discretization of time \citep{boland2017continuous}, and \textbf{(3)} a discretization of capacity.  Increasing the number of LA-neighbors of customers can tighten the LP relaxation by limiting the cycles (over customers localized in space) permitted in an LP solution, a principle also exploited by NG-routes \citep{baldacci2011new}. Discretization is parameterized by grouping capacity remaining (or time remaining into buckets).  For example given buckets associated with customer $u$ denoted $D_u$ we bin in the capacity remaining as follows:   $D_u=[[d_u,d_u+3],[d_u+4,d_u+5],[d_{u}+6,d_{u}+9]...[d_0-4,d_0] ]$.  Increasing the number of buckets makes the buckets smaller in range.  We refer to increasing the number of buckets as decreasing granularity while decreasing the number of buckets as increasing granularity.  Increasing the number of buckets in the time/capacity discretization sets tightens the LP relaxation by enforcing that routes must be time/capacity feasible more rigorously.  
However, increasing the number of LA-neighbors and decreasing the granularity of the time/capacity discretization increases the number of variables/constraints in the LP and hence increases the LP computation time.
Special selection of these sets (LA neighbors and time/capacity buckets) can ensure both fast LP optimization time, and a tighter LP relaxation relative to the baseline in \eqref{orig_opt}; and thus fast MILP optimization time. Given these sets (collectively called a parameterization) we need only solve the corresponding optimization problem as a MILP using an off-the-shelf solver to obtain an optimal solution to the VRPTW. 

We organize this section as follows. In Sections \ref{subsec_la_incor} and \ref{subsec_demand_incor}, we tighten the LP relaxation in \eqref{orig_opt} using LA-arcs, and time/capacity discretization respectively. In Section \ref{subsec_full_opt_incor}, we produce our novel MILP formulation of VRPTW. In Section \ref{subsec_sufficient}, we consider some properties of sufficient parameterizations, which we define as parameterizations that lead to the tightest bound possible given a maximum number of LA-neighbors per customer. In Section \ref{subsec_minimal}, we consider some properties of parsimonious parameterizations, which are parameterizations that have the widest time/capacity buckets and smallest LA-neighborhoods as possible, given fixed LP objective, leading to maximally efficient solutions to the LPs at each node of the branch \& bound tree of the corresponding MILP.  
\pagebreak
\subsection{Enforcing Elementarity via the Incorporation of Local Area Arcs}
\label{subsec_la_incor}
\begin{figure}
    \centering
    \includegraphics[width=0.95\textwidth]{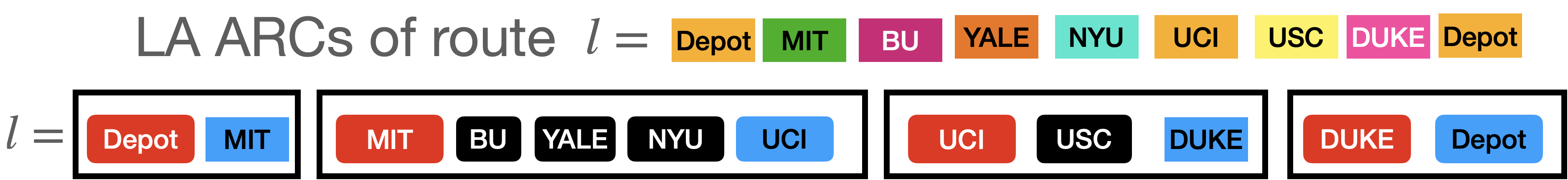}
    \caption{We display a route as a sequence of LA-arcs.  Each LA-arc describes a sequence of customers near the first customer in the arc, followed by a distant customer. The example uses American universities and LA-neighborhoods defined by university location. LA-arcs are color-coded as follows.  Red indicates the first customer (or depot) in an LA-arc, Black indicates intermediate customers in an LA-arc, and Blue indicates the final customer (or depot) in an LA-arc. The first customer of an LA-arc is shared with the last customer of the previous LA-arc.  Observe that the class of cycles permitted in a sequence of LA-arcs is highly limited because individual LA-arcs that make up the route must be elementary.  }
    \label{pic_LA_arc}
\end{figure}

In this section, we provide additional variables and associated constraints using Local Area (LA) arcs, which were introduced in the arXiv paper \citep{mandal2022local_2}.  LA-arcs are designed to eliminate cycles of customers localized in space as described by the $x$ variables used in \eqref{orig_opt}.
As a precursor, we provide some notation for describing LA-arcs, which we illustrate in Fig \ref{pic_LA_arc}.  For each $u\in N$ we use $N^{\odot}_u$ to denote the set of customers nearby $u$, not including $u$, called the LA-neighbors of $u$. Specifically $N^{\odot}_u$ is the set of $K$ closest customers to $u$ in space that are reachable from $u$ considering time and capacity; for user-defined $K$.\footnote{Other definitions can be used such as excluding from $N^{\odot}_u$ customers for which the route ordered as follows $\alpha,u,v,u,\bar{\alpha}$ is infeasible by time/capacity.} For ease of  notation we define $N^{\odot}_{+u}=\{N^{\odot}_{u}\cup u\}$.  We use $N^{\odot \rightarrow}_{u}$ to denote the subset of $N^+$ not in $N^{\odot}_{+u}$ or the starting depot meaning $N^{\odot \rightarrow}_{u} =N^+-(N^{\odot}_{+u}\cup \alpha)$.   For each $u\in N,v \in N^{\odot \rightarrow}_{u},\hat{N}\subseteq N^{\odot}_u$, we use $R^+_{u\hat{N}v}$ to denote the subset of all elementary sequences of customers (called orderings) starting at $u$; ending at $v$; including as intermediate customers $\hat{N}$; that are feasible with regard to time/capacity. Not all sets of intermediate customers $\hat{N}\subseteq N_u^{\odot}$ are feasible given $u,v$ due to time/capacity constraints.  For any given $u,v$ and set of intermediate customers in $\hat{N} \subseteq N_u$, we define $R_{u\hat{N}v}\subseteq R^+_{u\hat{N}v}$ to contain only those orderings on the efficient frontier with respect to \textbf{(a)} travel distance, \textbf{(b)} latest feasible departure time from $u$ and \textbf{(c)} earliest feasible departure time from $v$.  Feasibility of \textbf{(b)},\textbf{(c)} is a consequence of the time windows for the customers. Thus, orderings are preferred  with \textbf{(a)} lower cost, \textbf{(b)} that can be started later, or \textbf{(c)} that can be finished earlier. In practice, $R_{u\hat{N}v}$ is a small subset of all possible orderings $R^+_{u\hat{N}v}$ as most do not lie on the efficient frontier. 
We use $R_{u}$ to denote the union of $R_{u\hat{N}v}$ over all $v,\hat{N}$ where we clip off the final term in the ordering (meaning $v$). Observe that any optimal integer solution to VRPTW must select a route containing an ordering in $R_u$ followed immediately by some $v \in N^{\odot \rightarrow}_u$ for each $u \in N$. In practice, $R_{u}$ is smaller than the sum of the sizes of the $R_{u\hat{N}v}$ terms as many terms are replicated between $R_{u\hat{N}v}$ sets.  We provide a dynamic programming approach for generating $R_u$ in Section \ref{my_r_comp}.  

We describe an ordering $r\in R_u$ as follows.  For any $w \in N^{\odot}_{+u}$, $v \in N^{\odot}_u$ we define $a_{wvr}=1$ if $w$ immediately precedes $v$ in ordering $r$, and otherwise define $a_{wvr}=0$.  For any $w \in N^{\odot}_{+u}$ we define $a_{w*r}=1$ if $w$ is the final customer in ordering $r$, and otherwise define $a_{w*r}=0$. 

We use decision variables $y_r \in \mathbb{R}_{0+}\; \forall r \in R_u,\; u \in N$ to describe a selection of orderings to use in our solution.  We set $y_r=1$ for $r \in R_{u}$ to indicate that after servicing customer $u$ the sequence of customers described by ordering $r$ is used after, which a member of $N_u^{\odot \rightarrow}$ follows immediately, and otherwise set $y_r=0$.
We use $E^{*\odot}_{u}\subseteq E^*$ to denote the set of edges that can be contained in an ordering starting at $u$; thus $E^{*\odot}_{u}=\{ wv\in E^*, \mbox{ s.t. } w \in N^{\odot}_{+u},v\in N^{\odot}_u-w\} $.  We use $E^{*\odot}_{uw}$ to denote the set of edges in $E^*$ that can succeed the final customer $w$ in an ordering starting at $u$; thus $E^{*\odot}_{uw}=\{ wv \in E^* ; w \in N^{\odot}_{+u}, v \in N^{\odot\rightarrow}_u\}$.

We now add the following equations to \eqref{orig_opt}, which adds minimization over $y$ to enforce that the solution to the $x$ variables is consistent with the solution $y$, for which we provide an exposition below the equations.  
\begin{subequations}
\label{con2Group}
\begin{align}
 & \sum_{r \in R_{u}}y_{r}=1  \quad &\forall u \in N \label{oneLARCPer}\\   
  &  x_{wv}\geq \sum_{r \in R_u}a_{w,v,r}y_r \quad & \forall u \in N, wv \in E^{\odot *}_u
  \label{A_state_dir_local_con_1} \\
  &  
  \sum_{\substack{wv \in E^{*\odot}_{uw}}}
  x_{wv}\geq \sum_{r \in R_u}a_{w*r}y_r \quad &\forall u \in N,w \in N^{\odot}_{+u} \label{A_state_dir_global_con_1}
    \end{align}
    \end{subequations}
In \eqref{oneLARCPer} we enforce that one ordering in $R_u$ is selected to describe the activities succeeding $u$ for each $u \in N$.  
In \eqref{A_state_dir_local_con_1} we enforce that if an ordering $r \in R_u$ is selected in which $w$ is immediately followed by $v$; then it must be the case that $x_{wv}$ is selected meaning that $w$ is succeeded by $v$ in the solution over $x$.  In \eqref{A_state_dir_global_con_1} we enforce that if an ordering $r \in R_u$ is selected in which $w$ is the final customer in the ordering then, $w$ must be followed by a customer (or depot) that is not an LA-neighbor of $u$ (and not $u$) as described by $x$.  Observe that if $x$ is binary then $y$ must be binary as well.  The terms $y_{r}$ are not present in the objective for optimization.  The use of LA-arcs to remove cycles, is a key contribution as it allows the user to grow the size of the formulation as the computational capabilities of MILP solvers improve over time.  

It is possible to show (see Section \ref{la_parsimony}) that the number of LA-neighbors required can be less than that of $N^{\odot}_u$ for some $u$. To this end, we introduce $N_u\subseteq N_{u}^{\odot}$ to denote this reduced set.  We use $N_{+u}$ to denote $N_{u}\cup u$.  We use $N^{\rightarrow}_u$ to denoted the subset of $N^+$ that does not lie in $N_{+u}\cup \alpha$.  We use $E^{*}_{u}\subseteq E^*$ to denote the set of edges that can be contained in an ordering starting at $u$, thus $E^{*}_{u}=\{ wv\in E^*, \mbox{ s.t. }  w \in N_{+u},v\in N_u-w,\} $.  We use $E^{*}_{uw}$ to denote the set of edges in $E^*$ that can succeed the final customer $w$ in an ordering starting at $u$, thus $E^{*}_{uw}=\{ wv \in E^* ; v \in  N^{\rightarrow}_u\}$.
We write \eqref{con2Group} using LA-neighbors $N_u$ below.
\begin{subequations}
\label{con2Group1}
\begin{align}
 & \sum_{r \in R_{u}}y_{r}=1  \quad & \forall u \in N \label{oneLARCPer1}\\   
 &   x_{wv}\geq \sum_{r \in R_u}a_{w,v,r}y_r \quad & \forall u \in N,wv \in E^*_u \label{one_stateCon_11}\\
  &  
  \sum_{\substack{wv \in E^{*}_{uw}}}
  x_{wv}\geq \sum_{r \in R_u}a_{w*r}y_r \quad & \forall u \in N,w \in N_{+u} \label{state_dir_global_con_11}
    \end{align}
    \end{subequations}
Note that  $R_u$ is determined by $N^{\odot}_u$ not $N_u$; while the set of constraints enforced in \eqref{one_stateCon_11}-\eqref{state_dir_global_con_11} is defined by $N_u$.  Observe that this can create redundant $y$ terms. These are simply removed before optimization.  
\subsection{Capacity/Time Discretization Discovery}
\label{subsec_demand_incor}
This section considers using capacity/time discretization to tighten the LP relaxation in \eqref{orig_opt}. As mentioned earlier, this discretization is inspired by the work of \citep{boland2017continuous, boland2019price}. Earlier discretization methods can be found in \citep{appelgren1969column, appelgren1971integer, wang2002local}. 

Consider that we have a partition of capacity associated with each customer $u$ denoted $D_u$ into continuous ``buckets," e.g. $D_u=[[d_u,d_u+3],[d_u+4,d_u+5],[d_{u}+6,d_{u}+9]...[d_0-4,d_0] ]$. We index the buckets in  $D_{u}$ from smallest to greatest remaining demand with $k$; where $k=1$ is the bucket that contains $d_u$ and $k=|D_u|$ contains $d_{0}$. The lower/upper bound values of the $k'th$ bucket for customer $u$ are referred to as $d^-_{(u,k)},d^+_{(u,k)}$ respectively.

Consider a directed unweighted graph $G^D$ (which we also use for node-set) with edge set $E^D$.  
For each $u \in N$, $k\in D_u$ we create one node $i=(u,k)$ for which we alternatively write as $(u_i,d^-_i,d^+_i)$ where $u_i\leftarrow u, d^-_i\leftarrow d^-_k,d^+_i \leftarrow d^+_k$.  There is a single node  $(\alpha,d_0,d_0)$ also denoted $(\alpha,1)$ associated with the starting depot and a single node $(\bar{\alpha},0,d_0)$ also denoted $(\bar{\alpha},1)$ associated with the ending depot.  We now define the edges $E^D$.  
\begin{itemize}
\item 
For each $u$ we create a directed edge from $(\alpha,1)$ to $(u,|D_u|)$. 
 Traversing this edge indicates that a route starts with $u$ as its first customer.  
\item
For each $u$ we create a directed edge from $(u,1)$ to $(\bar{\alpha},1)$. Traversing this edge indicates that a route ends after having $u$ as its final customer.

\item For each $i=(u,k)$ for $k>1$ we connect $i$ to $j=(u,k-1)$.  Traversing this edge indicates that the vehicle servicing $u$ will not use all possible capacity.  These edges can be thought of as dumping extra capacity.  
\item For each pair of nodes $i=(u,k),j=(v,m)$ we connect $i$ to $j$ if the following are satisfied: $u\neq v$, $d^+_i-d_u\geq d^-_j$ and if $k>1$ we also require that $d^-_j> d^+_{(u,k-1)}-d_{u}$. Traversing this arc indicates that a vehicle has between $d^-_i$ and $d^+_i$ units of capacity remaining before servicing $u$ and, upon servicing, $u$ travels directly to $v$ where it has between $d^-_j$ and $d^+_j$ units of capacity remaining before servicing $v$. The second condition $d^-_j> d^+_{(u,k-1)}-d_{u}$ ensures that redundant edges are not generated.
\end{itemize}

We now describe flow formulation on the discretized capacity. This flow governs the amount of the resource capacity available at a given node. The following constraints are written below using decision variable $z^D_{ij}\in \mathbb{R}_{0+}$ to describe the edges traversed .
\begin{subequations}
\label{con_demand_group}
    \begin{align}
        \sum_{ij \in E^D}z^D_{ij}=\sum_{ji \in E^D}z^D_{ji} & \quad \forall i \in G^D; i=(u,k); u\notin [\alpha,\bar{\alpha}] \label{Con1A}\\
        x_{uv}=\sum_{\substack{ij \in E^D\\ i=(u,k)\\j=(v,m)}}z^D_{ij} & \quad \forall uv \in E^* \label{Con2A}
    \end{align}
\end{subequations}
In \eqref{Con1A}, we enforce that the solution $z^D$ is feasible with regard to capacity by enforcing that the flow incoming to a node is equal to the flow outgoing.    
In \eqref{Con2A}, we enforce that the solution $x$ is consistent with the solution capacity graph $z^D$. 

We now describe flow formulation on discretized time as we previously did for capacity.  Time buckets and associated graph and edge sets are denoted $T_u, G^T, E^T$ respectively. The constraints are defined analogously to capacity (with minimum/maximum values of $t^-_u$ and $t^+_u$ respectively in place of $d_u$ and $d_0$).  We write the corresponding constraints below.  We use non-negative decision variables $z^T_{ij}$ to describe the movement of vehicles. We set  $z^T_{ij}=1$ to indicate  that a vehicle leaves $u_i$ with time remaining in between $t^-_{i}$ and $t^+_i$ and leaves at $u_j$ with time remaining between $t^-_j$ and $t^+_j$.  We write the constraints for time analogous to \eqref{Con1A} and \eqref{Con2A} below.  

\begin{subequations}
\label{con_time_group}
    \begin{align}
        \sum_{ij \in E^T}z^T_{ij}=\sum_{ji \in E^T}z^T_{ji} & \quad \forall i \in G^T; i=(u,k); u\notin [\alpha,\bar{\alpha}] \label{Con1}\\
        x_{uv}=\sum_{\substack{ij \in E^T\\ i=(u,k)\\j=(v,m)}} z^T_{ij} &\quad \forall uv \in E^* \label{Con2}
    \end{align}
\end{subequations}
We refer to the finest possible granularity of capacity and time buckets for each $u\in N$ as  $D^{\odot}_u, T^{\odot}_u$.  
\subsection{Complete Linear Programming Relaxation}
\label{subsec_full_opt_incor}

Given a set of LA-neighbors, time buckets, and capacity buckets for each $u \in N$ denoted $N_u, T_u, D_u$ respectively, we write the MILP enforcing no-localized cycles in space (via \eqref{con2Group1}), capacity bucket feasibility (via \eqref{con_demand_group}), and time bucket feasibility (via \eqref{con_time_group}) as $\bar{\Psi}(\{ N_u, T_u, D_u \},\forall u \in N)$ below, with its LP relaxation denoted $\Psi(\{ N_u, T_u, D_u \},\forall u \in N)$.\\
\begin{subequations}
\label{comp_opt}
\begin{align}
       &\bar{\Psi}(\{ N_u ,T_u, D_u \},\forall u \in N)= \min_{\substack{x \in \{0,1\} \\ y\geq 0\\ z\geq 0\\t^+_u\geq \tau_u\geq t^-_u\\ d_0\geq  \delta_u \geq d_u}} 
\sum_{\substack{u \in N^+\\ v \in N^+}}c_{uv}x_{uv} \label{comp_obj_1} &\quad {\color{red} \mbox{Original Compact}}\\
       & \sum_{\substack{ uv \in E^*}}x_{uv}= 1 \quad \forall u \in N  \quad \label{serv_eq_1} &\quad {\color{red} \mbox{Original Compact}}\\
        & \sum_{vu \in E^*}x_{vu}=1  \quad  \forall u \in N \label{flow_con_1}& \quad {\color{red} \mbox{Original Compact}}\\
        & \delta_v-d_v \geq \delta_u-(d_0+d_v)(1-x_{vu}) \quad \forall u \in N , vu \in E^*  \label{demand_con_1}& \quad {\color{red} \mbox{Original Compact}}\\
        & \tau_v-t_{vu} \geq \tau_u-(t^+_u+t_{vu})(1-x_{vu}) \quad \forall u \in N, vu \in E^* \label{time_con_1} &\quad {\color{red} \mbox{Original Compact}}\\
                & \sum_{u \in N}x_{\alpha u}\geq \lceil \frac{\sum_{u \in N}d_u}{d_0}\rceil\label{minVeh_con2} &\quad  {\color{red} \mbox{Original Compact}}\\
    & \sum_{r \in R_{u}}y_{r}=1 \quad \forall u \in N &\quad \label{one_stateCon_1} {\color{blue} \mbox{LA-arc Movement Consistency}}\\
    & x_{wv}\geq \sum_{r \in R_u}a_{w,v,r}y_r \quad \forall u \in N,wv \in E^*_u
    &\quad  \label{state_dir_local_con_1}  {\color{blue} \mbox{LA-arc Movement Consistency}}\\
    &\sum_{\substack{wv \in E^*_{uw}}} x_{wv}\geq \sum_{r \in R_u}a_{w*r}y_r \quad \forall u \in N,w \in N_{+u} &\quad \label{state_dir_global_con_1} {\color{blue} \mbox{LA-arc Movement Consistency}}\\
    &\sum_{ij \in E^D}z^D_{ij}=\sum_{ji \in E^D}z^D_{ji} \quad \forall i \in G^D; i=(u,k); u\notin [\alpha,\bar{\alpha}] &\quad  {\color{cyan} \mbox{Flow Graph Capacity}} \label{frlow_grdem}\\
        &x_{uv}=\sum_{\substack{ij \in E^D\\i=(u,k)\\j=(v,m)}}z^D_{ij} \quad \forall uv \in E^* &{\color{cyan} \quad \mbox{Flow Graph Capacity}} \label{frlow_grdem2}\\
        &\sum_{ij \in E^T}z^T_{ij}=\sum_{ji \in E^T}z^T_{ji} \quad \forall i \in G^T; i=(u,k); u\notin [\alpha,\bar{\alpha}] & \quad  {\color{green} \mbox{Flow Graph Time}} \label{frlow_grtime_1}\\
        &x_{uv}=\sum_{\substack{ij \in E^T\\ i=(u,k)\\j=(v,m)}}z^T_{ji} \quad \forall uv \in E^* & \quad {\color{green} \mbox{Flow Graph Time}} \label{frlow_grtime_2}
    \end{align}
    \end{subequations}

To be able to use an off-the-shelf MILP solver to solve for \eqref{comp_opt} efficiently, we must have parameterization $(\{N_u,T_u,D_u\} \forall u \in N)$ that is as tight as possible and small as possible in terms of using fewer variables and constraints.  Since the tightest LP expressible is $\Psi(\{N^{\odot}_u,T^{\odot}_u,D^{\odot}_u\} \forall u \in N)$ we want $\Psi(\{N_u,T_u,D_u\} \forall u \in N)$ to be equal to $\Psi(\{N^{\odot}_u,T^{\odot}_u,D^{\odot}_u\} \forall u \in N)$ and refer to such parameterizations as sufficient. 

We want the LP $\Psi(\{N_u,T_u,D_u\} \forall u \in N)$ and hence each node in the MILP $\bar{\Psi}(\{N_u,T_u,D_u\} \forall u \in N)$ to be easily solved. One key condition for such parameterizations is that of being parsimonious; any further contraction of the sets in the parameterization decreases the LP value of the optimal solution. Another way to say this is that the parameterization has the widest buckets and smallest LA-neighborhoods possible, given fixed LP tightness, leading to maximally efficient solutions to the LPs at each node of the branch \& bound tree of the corresponding MILP. 
We refer to LPs that contain these two properties that as sufficient and parsimonious parameterization (SPP) and in Section \ref{sec_expand_contract} consider the construction of such LPs.

\subsection{Properties of Sufficient Parameterizations}
\label{subsec_sufficient}

We now consider some sufficient (but not always necessary) conditions of sufficient parameterizations. To this end we consider a candidate optimal LP solution $(x^*,y^*,z^*)$ given parameterization $(\{N^*_u,T^*_u,D^*_u\} \forall u \in N)$.

\paragraph{The property of bucket feasibility} applies to capacity and time.  This provides a condition for the LP solution that, if satisfied, guarantees that using the minimum granularity (minimum bucket sizes) of capacity would not tighten the bound given the current LA-neighbors. We require bucket feasibility to hold in our conditions for sufficient parameterizations.  
Bucket feasibility states no relaxation of capacity/time occurs as a consequence of using the current bucket parameters.  We formally define bucket feasibility for capacity and time below.
\begin{subequations}
\label{capTimeBucketFeas}
    \begin{align}
        &\textbf{Capacity Bucket Feasibility }[z^{D*}_{ij}>0]\rightarrow [d^+_i-d_{u_i}=d^+_j] \quad  \forall (i,j) \in E^D \; u_j  \notin \{ u_i,\bar{\alpha}\} \label{capBucketF}\\
       & \textbf{Time  Bucket Feasibility }[z^{T*}_{ij}>0]\rightarrow [\min(t^+_i-t_{u_i,u_j},t^+_j)=t^+_j] \quad \forall (i,j) \in E^T; u_j \notin \{ u_i,\bar{\alpha}\} \label{timeBucketF}
    \end{align}
\end{subequations}
Observe that if \eqref{capTimeBucketFeas} is satisfied, then no relaxation of capacity/time occurs as a consequence of using the buckets, hence using the finest capacity/time discretization $\{D^{\odot},T^{\odot}\}$ would not tighten the bound given the current LA-neighbors. We prove this formally in Appendix \ref{proofCapBucketFeas}.
%
%
%

\paragraph{The property of LA-arc feasibility} is a condition for the LP solution that, if satisfied, guarantees that using the $N_u^{\odot}$ for the LA-neighbors of each customer would not tighten the bound given the current time/capacity buckets.  Specifically, LA-arc feasibility states that given fixed $x^*$, a setting for $y$ obeys \eqref{con2Group}, not merely \eqref{con2Group1}. We require LA-arc feasibility to hold in our conditions for  sufficient parameterizations.
\subsection{Properties of Parsimonious Parameterizatons}
\label{subsec_minimal}

We now consider some necessary but not always sufficient properties of parsimonious parameterizations.  

\subsubsection{Capacity/Time Parsimony}
\label{demMinimal}
Consider the dual solution to the LP in \eqref{comp_opt}, which we denote $\pi$ where $\pi^D_i$ is the dual value associated with \eqref{frlow_grdem} for a given $i$.  We write the reduced cost of $z^{D}_{ij}$ as $\bar{z}^D_{ij}$ for any $i,j$ s.t. $u_i=u_j$.  Observe that $\bar{z}^D_{ij}=-\pi^D_i+\pi^D_j$.  

Now observe that merging capacity buckets $i,j$ corresponds to enforcing that $\pi^D_i=\pi^D_j$. Here merging $i=(u_i,k+1)$ with $j=(u_i,k)$ corresponds to removing buckets corresponding to $(u_i,k+1)$ and $(u_i,k)$ and replacing them with a bucket associated with $u_{i}$ and range of capacity $[d^-_{(u_i,k)},d^+_{(u_i,k+1)}]$.  When $\pi^D_i=\pi^D_j $ for $u_i=u_j$ and $i,j \in E^D$ is satisfied, then observe that merging capacity buckets associated with $i,j$ together does not weaken the LP relaxation in \eqref{comp_opt}. However, it does decrease the number of constraints/variables in \eqref{comp_opt}, leading to faster optimization of the LP.

This same logic holds for time buckets.  Let $\pi^T_{i}$ denote the dual variable associated with  \eqref{frlow_grtime_1}.  We write the reduced cost for variable $z^T_{ij}$, which we denote as  $\bar{z}^T_{ij}$ as follows $\bar{z}^T_{ij}=-\pi^T_i+\pi^T_j$ for any $z^T_{i,j}$ s.t. $u_i=u_j$ and $i,j \in E^T$.  Thus, when $\pi^T_i=\pi^T_j$ for $u_i=u_j$ observe that we can merge the time buckets $i,j$ without loosening the bound in \eqref{comp_opt}.  However, as is the case for capacity buckets, this decreases the number of constraints/variables in \eqref{comp_opt}, leading to faster optimization of the LP.   
 
\subsubsection{LA Neighborhood Parsimony}
\label{la_parsimony}
Let us now consider the construction of parsimonious LA-neighborhoods to enforce in \eqref{comp_opt}. Let us consider the following additional notation.
Let $N_u^{k}$ be the $k$ closest LA-neighbors in $N_u$ to $u$ (in terms of distance).  We use $ N_{+u}^{k}$ to denote $N_u^{k}\cup u$.   Here $k$ ranges from $0$ to $|N_u|$.  

For any ordering $r$ where $N_r$ are the customers in the order $r$ and $w \in N_r$ we define $a^k_{wr}=1$ if all customer in $r$ prior to and including $w$ lie in $N_{+u}^k$ and otherwise set $a^k_{wr}=0$.

For each $u \in N,r \in R_u,w\in N_{+u}^k,v\in N_u^k-w,k$ we define $a_{wvr}^k=1$ if $a_{wvr}=1$ and $a^k_{vr}=1$.  
For each $u,r \in R_u,w\in N_{+u}^k$ we define $a_{w*r}^k=1$ if $a^k_{wr}=1$ and either ($a_{w*r}=1$ or $\sum_{v \notin N_u^k-w }a_{wvr}=1$ ) ; and otherwise set $a_{w*r}^k=0$.

We now replace \eqref{state_dir_local_con_1} and \eqref{state_dir_global_con_1} with the following inequalities parameterized by tiny positive value $\epsilon$ and dual variables noted in brackets. We use slack constant $\epsilon$  which we multiply by the number of LA-neighbors associated with a given constraint in order to encourage the LP solver to use constraints corresponding to as few LA-neighbors as possible. Later this facilitates the construction of a more parsimonious LP as larger LA neighborhoods are not be used in the generated LP if the dual variables are not active.    

\begin{subequations}
\label{conGroupK}
    \begin{align}
     &       \epsilon k+x_{wv}\geq \sum_{r \in R_u}a^k_{w,v,r}y_r \quad \forall u \in N,w \in N^k_{+u},k \in [1,2,...|N_u|],v\in N^k_u,wv \in E^* \label{state_dir_local_con_2} \quad [-\pi_{uwvk}]  \\
   & \epsilon k+\sum_{\substack{v \notin N^+-(N^k_{+u})\\ wv \in E^*}}x_{wv}\geq \sum_{r \in R_u}a^k_{w*r}y_r \quad \forall u \in N,k \in [1,2,...|N_u|],w \in N^k_{+u}\label{state_dir_global_con_2} \quad [-\pi_{uw*k}]  
    \end{align}
\end{subequations}

We now solve the LP form in \eqref{comp_opt} replacing \eqref{state_dir_local_con_1} and \eqref{state_dir_global_con_1} with \eqref{state_dir_local_con_2} and \eqref{state_dir_global_con_2} respectively.  We refer to the associated LP as $\Psi^*(\{N_u,T_u,D_u\} \forall u \in N)$

Observe that the dual LP for $\Psi^*(\{N_u,T_u,D_u\} \forall u \in N)$ has a penalty for using dual variables of primal constraints associated with larger LA-neighborhoods than necessary.  Consider any $u$, and let $k_u$  be the largest neighborhood size containing a dual variable of the form  \eqref{state_dir_local_con_2} or \eqref{state_dir_global_con_2} with positive value. We write $k_u$ formally below. 
\begin{align}
\label{applyContract_LA}
k_u\leftarrow \max_{|N_u|\geq k\geq 1 } k[0<( (\sum_{w \in N_u}\pi_{uwk})+ (\sum_{\substack{w \in N^k_{+u}\\ v\in N^k_u-w\\ wv \in E^*}}\pi_{uwvk}))] \quad \forall u \in N
\end{align}
If $k_u<|N_u|$, then we can decrease the size of $N_u$ to $k_u$ without loosening the bound.  This decreases the number of constraints in \eqref{comp_opt} of the forms \eqref{state_dir_local_con_1}, and \eqref{state_dir_global_con_1}.  This also decreases the number of variables that need to be considered as some of the $y_r$ variables become redundant. This, in turn, leads to the faster resolution of the associated LP. Observe that given $u$, if no dual variables of the form \eqref{conGroupK} are positively valued, then setting $N_u$ to be empty does not loosen the bound. 

\section{LA-Discretization: An Algorithm for a Sufficient and Parsimonious Parameterization Construction}
\label{sec_expand_contract}

In this section, we consider the construction of a sufficient and parsimonious parameterization (SPP). We construct an SPP by iterating between the following two steps:  \textbf{(1)} Solving $\Psi^*(\{N_u,T_u,D_u\} \forall u \in N)$ and \textbf{(2)} Augmenting $\Psi^*(\{N_u,T_u,D_u\}$ to achieve sufficiency.  Whenever \textbf{(1)} tightens the bound (with respect to the previous iteration), we attempt to achieve parsimony by decreasing the size of some terms in $\{N_u,T_u,D_u\}$ for some customers. This decrease is always done not to loosen the bound but to achieve the specific necessary conditions for parsimony.  We only aim towards parsimony after improving the bound (not at every step) to avoid cycling in the LP.  

Thresholds determine the buckets, which we denote as $W^D_u,W^T_u$ (sorted from smallest to largest) for capacity and time for customer $u$, respectively, where the thresholds determine the maximum value of the buckets. For example if the thresholds for $W^D_u$ are $[d_u+5,d_u+12]$ then the associated buckets are $[d_u,d_u+5],[d_u+6,d_u+12],[d_u+13,d_0]$.

We organize this section as follows. In Section \ref{sec_alg_contract_time_demand}, we consider the contraction and expansion operations for time/capacity buckets, while Section \ref{sec_alg_contract_LA} does the same for the number of LA-neighbors. In Section \ref{sec_alg_complete}, we describe our complete optimization algorithm, which we call LA-Discretization.
 \subsection{Contracting/Expanding Capacity/Time Buckets}
\label{sec_alg_contract_time_demand}
Recall from Section \ref{demMinimal} that for any $i,j \in E^D$ for which $u_i=u_j$ if $\pi^D_i=\pi^D_j$ then we can merge nodes $i,j$ without loosening the bound and remove $d^+_j$ from $W^D_{u_j}$.  The same logic applies to time. For any $i,j \in E^T$ for which $u_i=u_j$ if $\pi^T_i=\pi^T_j$, we can merge nodes $i,j$ without loosening the bound and remove $t^+_j$ from $W^T_{u_j}$. 

To enforce  bucket feasibility with respect to capacity we do the following.  For each $i,j$ s.t. $z^D_{ij}>0$ where $i=(u_i,d^-_i,d^+_i)$ and $j=(u_j,d^-_j,d^+_j)$ and $j\neq \bar{\alpha}$ and $u_i\neq u_j$: 
 We add to $W^D_{u_j}$ the term $d^+_i-d_{u_i}$ (if this term is not already present in $W^D_{u_j}$ ).  Observe that this operation can never loosen the relaxation as undoing it simply corresponds to requiring that  $\pi^D_{i}=\pi^D_j$ in the dual.  

The same approach applies to enforcing bucket feasibility with respect to time.  
For each $i,j$ s.t. $z^D_{ij}>0$ where $i=(u_i,t^-_i,t^+_i)$ and $j=(u_j,t^-_j,t^+_j)$ and $j\neq \bar{\alpha}$:  We add to $W^T_{u_j}$ the term $\min(t^+_i-t_{u_i,u_j},t^+_j)$ (if this term is not already present in $W^T_{u_j}$ ). Observe that this operation can never loosen the relaxation as undoing it simply corresponds to requiring that  $\pi^T_{i}=\pi^T_j$ in the dual.
\subsection{Contracting/Expanding Local Area Neighborhoods}
\label{sec_alg_contract_LA}
Recall from  \eqref{applyContract_LA} that we can set the number of LA-neighbors to the smallest number for which there is an associated non-zero dual value without loosening the bound.  
Thus our contraction operation on LA-neighborhoods is written as follows.  
\begin{align}
\label{contract_LA_opt}
   |N_u|\leftarrow  \max_{|N_u|\geq k\geq 1 } k[0<( (\sum_{w \in N_u}\pi_{uwk})+ (\sum_{\substack{w \in N^k_{+u} \\ v\in N^k_u-w\\wv \in E^*}}\pi_{uwvk}))] \quad \forall u \in N
\end{align}
To expand the neighborhoods, we can iterate over $u\in N$, then, for each $u$, identify the smallest number of LA-neighbors required to induce a violation of the current LP. To determine if a given number of LA-neighbors $k$ induces a violation of the current LP we simply check to see if a fractional solution $y_r$ exists to satisfy \eqref{state_dir_global_con_1}, \eqref{one_stateCon_1} and \eqref{state_dir_local_con_1} where we fix the $x$ terms to their current values and $N_u$ to $N^k_u$.  

A simpler approach is as follows.  We can periodically simply reset the number of LA-neighbors to the maximum number $|N^{\odot}_u|$ for each $u\in N$.  Next, we solve the LP and apply a contraction operation described in \eqref{contract_LA_opt}. Our experiments employ this approach when expansion operations on capacity/time buckets have not recently tightened the LP bound.  


\subsection{Complete Algorithm: LA-Discretization}
\label{sec_alg_complete}

We now consider the construction of an SPP by iteratively solving the LP relaxation and tightening the LP relaxation to aim toward sufficiency; each time we tighten the relaxation, we decrease the size of the parameterization with parsimony. We provide our algorithm in Alg \ref{final_opt} with exposition below.
\begin{algorithm}[!b]
 \caption{LA-Discretization Algorithm}
\begin{algorithmic}[1] 
\State $\{N_u^{\odot},N_u,W^D_u,W^T_u \quad \forall u \in N\} \leftarrow $ from User \label{init1A}
\State iter\_since\_reset$\leftarrow$ 0 \label{init_state_iter}
\State last\_lp\_val$\leftarrow -\infty$ \label{initLP}
\Repeat \label{startRepeat}
\If {iter\_since\_reset $\geq \zeta$}: \label{reset_start}
    \State $N_u \leftarrow N^{\odot}_u \quad \forall u \in N$ 
\EndIf \label{restart_finish}
\State $[z,y,x,\pi,$lp\_objective$]\leftarrow$ Solve $\Psi^*(\{ N_u,T_u,D_u \} \forall u \in N)$ \label{startLP}
\If {lp\_objective$>$last\_lp\_val+MIN\_INC} 
\label{start_contraction}
\State   $W^{D}_u\leftarrow $ Merge any nodes $i,j$ for which $i,j \in E^D$,$u_i=u_j=u$ and $\pi^D_i=\pi^D_j$ for all $u \in N$ \label{start_cut_meat} 
\State   $W^{T}_u\leftarrow $ Merge any nodes $i,j$ for which $i,j \in E^T$,$u_i=u_j=u$ and $\pi^T_i=\pi^T_j$ for all $u \in N$ \label{start_cut_meat_1} 
\State $\{ N_u \quad \forall u \in N\} \leftarrow $ Apply \eqref{contract_LA_opt} to contract LA-neighbors.   \label{LACompress}
\State last\_lp\_val$\leftarrow$ lp\_objective\label{lp_step_update}
\State iter\_since\_reset $\leftarrow 0$ \label{reset_line}
\EndIf \label{end_contraction}
\State  $W^T_{u}\leftarrow W^T_{u} \cup_{ij \in E^T ;u_j=u; z^T_{ij}>0} \min(t^+_i-t_{u_i,u_j},t^+_j)$ for all $u\in N$\label{startFor11}
\State  $W^D_{u}\leftarrow W^D_{u} \cup_{ij \in E^D;u_j=u; z^D_{ij}>0} (d^+_i-d_{u_i})$ for all $u\in N$ \label{startFor00}

\State  iter\_since\_reset ++ \label{inc_inter_step}
\Until{$\{N_u,W^D_u,W^T_u\quad \forall u \in N\} $ are unchanged from the previous iteration OR iter\_since\_reset$>$ITER\_MAX} \label{endRepeart}
\State   $W^{T}_u,W^{D}_u,N_u\leftarrow $ via contraction operations from  Lines \ref{start_cut_meat}-\ref{LACompress} \label{start_cut_meat2}
\State $x\leftarrow $Solve $\bar{\Psi}(\{ N_u,T_u,D_u \} \forall u \in N)$ as MILP (only x is integer) \label{finOLP}
\end{algorithmic}
\label{final_opt}
\end{algorithm} 

\begin{itemize}
\item Line \ref{init1A}: Initialize the LA-neighbor sets, time buckets and capacity buckets. 
In our implementation We initialize $N_u \leftarrow N^{\odot}_u$ for all $u \in N$.  We initialize $W^D_u$ parameterized by a bucket size $d^s=5$ where $W^D_u=\{ [d_u,d_u+d^s), [d_u+d^s,d_u+2d^s)...)$. The last bucket may be smaller than $ d^s$ .  We initialize $W^T_u$ parameterized by a bucket size $t^s=50$ where $W^T_u=\{ [t^-_u,t^-_u+t^s), [t^-_u+t^s,t^-_u+2t^s)...)$. The last bucket may be smaller than $t^s$. We define the number of elements in $N^{\odot}_u$ for each $u \in N$ equal to user defined parameter $n^s$ where $n^s=6$ in our experiments. 
\item Line \ref{init_state_iter}-\ref{initLP}: We set the number of iterations since the last contraction operations, which we denote as iter\_since\_reset to zero.  We also set the incumbent LP value denoted  last\_lp\_val LP objective to $-\infty$. This value last\_lp\_val is the LP value immediately after the most recent contraction operation which is $-\infty$ since no contractions have yet happened.  
\item Line \ref{startRepeat}-\ref{endRepeart}: Construct a sufficient parameterization. 
\begin{itemize} 
\item Line \ref{reset_start}-Line \ref{restart_finish}:  If we have not tightened the bound by user defined parameter MIN\_INC (MIN\_INC=1 in experiments) for a given amount of iterations (denoted $\zeta$), then we aim towards sufficiency by increasing the LA-neighborhoods of all customers to their maximum size. We set $\zeta=9$ in our experiments.  
\item Line \ref{startLP}: Solve the LP given the current parameterization.   
\item Line \ref{start_contraction}-\ref{end_contraction}: 
If we have improved the LP relaxation by MIN\_INC since the last contraction operation, we apply contraction operations to aim toward parsimony for our parameterization.  
\begin{itemize}
\item Line \ref{start_cut_meat}-\ref{LACompress}: Apply contraction operations to time/capacity buckets and LA-neighborhoods.  
\item Line \ref{lp_step_update}: Update the LP objective immediately after the most recent contraction operation.
\item Line \ref{reset_line}:  Set the number of iterations since the most recent contraction operation to zero. 
\end{itemize}
\item Line \ref{startFor11},\ref{startFor00}: Aim towards sufficiency by applying expansion operations to capacity and time buckets, respectively.  
\item Line \ref{inc_inter_step}:  Increment iter\_since\_reset.
\item Line \ref{endRepeart}:  Terminate optimization when sufficiency is achieved, or progress towards sufficiency has ceased.  We stop the algorithm after ITER\_MAX(=10 in experiments) consecutive iterations of failing to tighten the LP relaxation by a minimum of MIN\_INC (=1  in experiments) combined. In Section \ref{sec_proof_converge} we prove that if ITER\_MAX=$\infty$ then Alg \ref{final_opt} (Lines \ref{startRepeat}-\ref{endRepeart}) produces a sufficient parameterization.
\end{itemize}
\item Lines \ref{start_cut_meat2}: Decrease the LP size by modifying the LP for parsimony.  
\item Line \ref{finOLP}: Finally solve the appropriate MILP in \eqref{comp_opt} given our SPP. Note that only $x$ is enforced to be binary, and all other variables are continuous. Enforcing $y$ to be binary does not change the MILP solution objective but may accelerate optimization depending on the specifics of the MILP solver. We enforced $y$ to be binary in experiments.   
\end{itemize}
%


\section{Local Area Arc Computation}
\label{my_r_comp}

In this section we consider the computation of the efficient frontier of orderings $R_{u\hat{N}v}$.  We define helper terms $P$ to denote the set of all unique values of $u \in N$,$\hat{N}\subseteq N^{\odot}_u$,$v \in N^{\odot \rightarrow}_u$, which we index by $p$.  Here we describe a specific $p$ as $(u_p,N_p,v_p)$. 
We organize this section as follows.  
In Section \ref{subsec_make_term_notation}, we provide notation needed to express the material in this section.  In Section \ref{subsec_make_term_recur}, we describe a recursive relationship for the cost of an LA-arc $p$ given a departure times from $u_p,v_p$.  In Section \ref{subsec_make_term_fronteir}, we introduce the concept of an efficient frontier of orderings for each $p\in P$ trading off latest feasible departure time from $u_p$ (later is preferred); earliest departure time from $u_p$ when no waiting is required (earlier is preferred); and cost (lower cost is preferred).  In Section \ref{subsec_make_term_alg}, we provide a dynamic programming based algorithm to compute the efficient frontier for all LA-arcs jointly.  In Section \ref{sec_proof}, we provide a proof referenced in Section \ref{subsec_make_term_recur}.
\subsection{Notation for LA-Arcs Modeling Time}
\label{subsec_make_term_notation}
Let $P^+$ be a super-set of $P$ defined as follows.  Here $p=(u_p,N_p,v_p)$ for $u_p \in N,N_p\subseteq N,v_p\in N^{+}$ lies in $P^+$ if and only if all of the following conditions hold. There exists a $u\in N$ s.t. $u_p \in (u \cup N_u);\; v_p\in N^{\odot \rightarrow }_u; \; N_p \subseteq N_u-(u_p+v_p)$.  We use helper term $c_{p,t_1,t_2}$ to denote the cost of lowest cost feasible ordering departing $u_p$ at time $t_1$, departing $v_p$ at $t_2$, visiting $N_p$ in between; where $r_{p,t_1,t_2}$ denotes the associated ordering. We define $R_{u\hat{N}v}$ as the union of $r \in R^+_{u\hat{N}v}$ for which $r=r_{p,t_1,t_2}$ for some $t_1,t_2$.  
We use $R^+_p,R_p$ as short hand for $R^+_{u_pN_pv_p}$ and $R_{u_pN_pv_p}$ respectively. We describe $r$ as an ordered sequence listed from in chronological order of visits as follows $[u^r_1,u^r_2,u^r_3,...v^r]$ with the final element denoted $v^r$. We define an ordering $r^-$, to be the same as $r$ except with $u^r_1$ removed, meaning $r^-=[u^r_2,u^r_3,...v^r]$.  

For any given time $t$ and ordering $r\in R^+_p$, we define $T_{r}(t)$ as the earliest time a vehicle could depart $v^r$, if that vehicle departs $u^r_1$ with $t$ time remaining and follows the ordering of $r$.  We write $T_r(t)$ mathematically below using $-\infty$ to indicate infeasibility, and $[]$ to denote the binary indicator function. 
\begin{align}
\label{trtDef}
    T_{r}(t)=T_{r^-}(-t_{uw}+\min(t,t^+_u))-\infty[t<t^-_u] \quad \forall r\in R^+_p,p \in P^+, u^r_1=u,w=u^r_2
\end{align}
Given $T_{r}(t)$ we express $c_{p,t_1,t_2}$ as follows.
\begin{align}
\label{lookupTravelTime_orig0xx}
c_{p t_1 t_2}=\min_{\substack{r \in R^+_p\\ T_{r}(t_1)\geq t_2}}c_{r} \quad \forall p\in P^+,t_1,t_2
\end{align}
Here $c_{p t_1 t_2}=\infty $ if no $r$ exists satisfying $T_{r}(t_1)\geq t_2$.  Observe that  using \eqref{lookupTravelTime_orig0xx} to evaluate $c_{p t_1 t_2}$ is challenging as we would have to repeatedly evaluate $T_{r}(t)$ in a nested manner.  
\subsection{Recursive Definition of Departure Time}
\label{subsec_make_term_recur}
In this section we provide an alternative characterization of $T_r(t)$, that later provides us with a mechanism to efficiently evaluate $c_{p t_1 t_2}$.  For any $r\in R^+_p,p\in P^+$ we describe $T_{r}(t)$ using the helper terms $\phi_{r},\hat{\phi}_{r}$ defined as follows.   %
We define $\phi_{r}$ to be the earliest time that a vehicle could leave $u^r_{1}$ without waiting at any customer in $u_p \cup N_p\cup v_p$ if $t^-$ terms were ignored (meaning all $t^-$ terms are set to $-\infty$). Similarly we define $\hat{\phi}_{r}$ as the latest time a vehicle could leave $u^r_1$ if $t^+$ terms were ignored (meaning all $t^+$ terms are set to $\infty$).   
Below we define $\phi_{r},\hat{\phi}_{r}$ recursively.
\begin{subequations}
\label{tauDef}
\begin{align}
\phi_{r}=\min(t^+_u,t^+_v+t_{uv}) \quad \forall r \in R^+_p, p=(u,\{ \},v),p \in P^+\\
    \hat{\phi}_{r}=\max(t^-_u,t^-_v+t_{uv}) \quad \forall r \in R^+_p, p=(u,\{ \},v),p \in P^+\\
    \phi_{r}=\min(t^+_u,\phi_{r^-}+t_{uw}) \quad \forall r \in R^+_p, u^r_{1}=u,u^r_{2}=w, |N_p|>0, p \in P^+\\ 
    \hat{\phi}_{r}=\max(t^-_u,\hat{\phi}_{r^-}+t_{uw}) \quad \forall r \in R^+_p,u^r_{1}=u,u^r_{2}=w, |N_p|>0, p \in P^+ 
\end{align}
\end{subequations}
Let $c_r$ denote the total travel distance on ordering $r$.  We rewrite $T_{r}(t)$ using $\phi_{r},\hat{\phi}_{r}$, with intuition provided subsequently (with proof of equivalence in Appendix \ref{sec_proof}).   
\begin{align}
\label{ituDef}
    T_{r}(t)=-c_r+\min(t,\phi_{r}) -\infty[\min(t,t^+_{u_p})<\hat{\phi}_{r}] \quad \forall p \in P^+,r \in R^+_p 
\end{align}
We now provide an intuitive explanation of \eqref{ituDef}.  Observe that leaving $u^r_1$ after $\hat{\phi}_{r}$ and following the ordering $r$ is infeasible. 
 Observe that leaving $u^r_1$ at time $t$ and following ordering $r$ for $t> \phi_{r}$ incurs a waiting time of $t-\phi_{r}$ over the course of the ordering, and otherwise incurs a waiting time of zero. Thus we subtract the travel time $c_r$ from $\phi_{r}$ to obtain the departure time at $v^r$.  Given $R^+_p$ and \eqref{ituDef} we write $c_{p,t_1,t_2}$ below.  
\begin{align}
\label{lookupTravelTime_orig0}
c_{p t_1 t_2}=\min_{\substack{r \in R^+_p\\ t_1\geq \hat{\phi}_{r}\\ t_2\leq -c_r+\min(t_1,\phi_{r})}}c_{r} \quad \forall p\in P^+,t_1,t_2
\end{align}
If no $r$ satisfies the constraints in \eqref{lookupTravelTime_orig0} then $c_{p t_1 t_2}=\infty$.  Observe that $|R^+_p|$ can grow factorially with $|N_p|$ thus using \eqref{lookupTravelTime_orig0} is impractical when $|N_p|$ grows.
\subsection{Efficient Frontier}
\label{subsec_make_term_fronteir}
In this section we describe a sufficient subset of $R^+_p$ denoted $R_p$ s.t. that applying \eqref{lookupTravelTime_orig0} over that subset produces the same result as if all of $R^+_p$ were considered.  
Given $p$, a necessary criterion for a given  $r \in  R^+_p$ to lie in $ \in R_p$, is that $r$ is not Pareto dominated by any other such ordering in $R^+_p$ with regards to latest time to start the ordering (later is preferred), cost (lower is preferred), and earliest time to start the ordering without waiting (earlier is preferred).  Thus, we prefer smaller values of $\hat{\phi}_{r},c_r$ 
and larger values of $\phi_{r}$.  
We use $R_p\subseteq R^+_p$ to denote the efficient frontier of $R^+_p$ with regards to $\hat{\phi}_{r},c_r,\phi_{r}$. Here for any $r\in R^+_p$ we define $r$ to lie in $R_p$ IFF no $\hat{r}\in R^+_{p}$ exists for which all of the following inequalities hold, (and at least one holds strictly meaning not with an equality ).  
    \textbf{(1)}$\hat{\phi}_{r}\geq \hat{\phi}_{\hat{r}}$.  
    \textbf{(2)}$c_{r}\geq c_{\hat{r}}$.    
    \textbf{(3)}$\phi_{r}\leq \phi_{\hat{r}}$.  We use the following stricter criteria for \textbf{(3)}: there is no time $t$ when we leave $u_p$ using $r$ in which we could depart $v_p$ before we could depart $v_p$ using $\hat{r}$. This criteria is written as follows $\phi_{r}-c_{r}\leq \phi_{\hat{r}}-c_{\hat{r}}$.  Observe that given the efficient frontier $R_p$ that we can compute $c_{pt_1t_2}$ as follows.  
\begin{align}
\label{lookupTravelTime}
c_{p t_1 t_2}=\min_{\substack{r \in R_p\\ t_1\geq \hat{\phi}_{r}\\ t_2\leq -c_r+\min(t_1,\phi_{r})}}c_{r} \quad \forall p\in P^+,t_1,t_2
\end{align}
If no $r$ satisfies the constraints in \eqref{lookupTravelTime} then $c_{p t_1 t_2}=\infty$.
Observe that we can not construct $R_p$ by removing elements from $R^+_p$ as enumerating $R^+_p$ would be too time intensive. 
\subsection{An Algorithm to Generate the Efficient Frontier}
\label{subsec_make_term_alg}

In this section we provide an efficient algorithm to construct all $R_p$ jointly without explicitly enumerating $R^+_p$.  
In order to achieve this we exploit the following observation for any $p\in P^+$ s.t. $|N_p|>0$.  Observe that if $r\in R_p$ then $r^-$ must lie in $R_{\hat{p}}$ where $\hat{p}=(u_{\hat{p}}=u^r_2,N_{\hat{p}}=N_{p}-u^r_2,v_{\hat{p}}=v^r) $. 

We construct $R_p$ terms by iterating over $p\in P^+$ from smallest to largest  with regard to $|N_p|$, and constructing $R_p$ using the previously generated efficient frontiers.  For a given $p$ we first add $u_p$ in front of all possible $r^-$ (called predecessor orderings) 
; then remove all orderings 
that do not lie on the efficient frontier.  The set of predecessor orderings for a given $r$ is written below.  
\begin{align}
    \cup_{\substack{w \in N_p\\ \hat{p}=(w,N_p-w,v_p)}}R_{\hat{p}}
\end{align}
Observe that the base case where $N_p$ is empty then there is only one member of $R_p$ which is defined by sequence $r=[u_p,v_p]$ when feasible, and otherwise is empty. In Alg \ref{better_ver} we describe the construction of $R_p$, which we annotate below.
\begin{algorithm}
\caption{Computing the Efficient Frontier for all $p\in P$ }
\begin{algorithmic}[1]
\For{$p\in P^+$ from smallest to largest in terms of $|N_p|$ with $|N_p|>0$} \label{outer_0}
\State $R_{p}\leftarrow \{ \}$ \label{init_rp}
\For{$w \in N_p; \hat{p}=(w,N_p-w,v_p),r^- \in R_{\hat{p}}$} \label{outer1}
\State $r\leftarrow [u_p,r^-]$ \label{grab_gplus}
\State Compute $c_{r}$, and $\phi_{r}$, $\hat{\phi}_{r}$ via \eqref{tauDef} \label{comp_tau}
\If{$\hat{\phi}_{r} \le t^+_u$} \label{if_add}
\State $R_{p}\leftarrow R_p \cup r$  \label{augment_rp}
\EndIf \label{if_add_e}
\EndFor \label{outer1_e}
\For{$r \in R_p,\hat{r} \in R_p-r$}\label{rem_dominated}
\If{$c_{r}\geq c_{\hat{r}}$ and $\hat{\phi}_{r}\geq \hat{\phi}_{\hat{r}}$ and $\phi_{r}-c_{r}\leq \phi_{\hat{r}}-c_{\hat{r}}$ and at least one inequality is strict} \label{if_0}
\State $R_{p}\leftarrow R_{p}-r $ \label{subtrTerm} 
\EndIf \label{if_e}
\EndFor \label{rem_dominated_e}
\EndFor \label{outer_0_e}
\end{algorithmic}
\label{better_ver}
\end{algorithm}
\begin{itemize}
    \item Line \ref{outer_0}-\ref{outer_0_e}:  Iterate over $p \in P^+$ from smallest to largest with regards to $|N_p|$ and given $p$ compute $R_p$.
    \begin{itemize}
    \item Line \ref{init_rp}: Initialize $R_p$ to be empty.
    \item Line \ref{outer1}-\ref{outer1_e}:  Compute all possible terms for $R_p$ by adding a customer $u_p$ in front of all potential predecessor $r^-$ terms.
    \begin{itemize}
        \item Line \ref{grab_gplus}:  Add $u_p$ to the front of $r^-$ creating $r$.
        \item Line \ref{comp_tau}:  Compute $c_{r}$, and  $\phi_{r}$, $\hat{\phi}_{r}$ via \eqref{tauDef}.
        \item Line \ref{if_add}-\ref{if_add_e}:  If $r$ is feasible for some departure time $t_1$ then we add $r$ to $R_p$.
    \end{itemize}
    \item Line \ref{rem_dominated}-\ref{rem_dominated_e}:  Iterate over members of $R_p$ and remove terms that do not lie on the efficient frontier.   
    \begin{itemize}
    \item Line \ref{if_0}-\ref{if_e}:  If $r$ is dominated by $\hat{r}$ then we remove $r$ from $R_p$.  
    \end{itemize}
\end{itemize}
\end{itemize}
We compute $R_p$ for all $p$ using Alg \ref{better_ver} once prior to running Alg \ref{final_opt}.  To construct $R_u$ we then take the union of all terms in $R_{u\hat{N}v}$ for all $v\in N^{\odot \rightarrow }_p, \hat{N}\subseteq N_u^{\odot}$ then finally clip $v$ from the end.

\section{Experiments}
\label{sec_exper}

In this section, we provide experimental validation for LA-Discretization. We compare LA-Discretization to the baseline optimization of solving the two-index compact MILP in \eqref{orig_opt}. We validated LA-Discretization on the standard Solomon instance data sets 
\citep{solomon1987algorithms} with additional results in Extension \ref{allResults}.
We provide the maximum computation time for the MILP call of 1000 seconds for each solver. All experiments are run using Gurobi with default parameters, and all codes are implemented in Python. We used a 2022 Apple Macbook Pro with 24 GB of memory and an Apple M2 chip.  As is often done in the literature, all distances are rounded down to the nearest tenth.

We fixed the optimization parameters to reasonable values fixed once. These values are shared by all problem instances on all data sets.  We set $n^s=6, d^s=5, t^s=50$,  ITER\_MAX=10, $\zeta=9$ and MIN\_INC=1. We did not find the solver to be sensitive to these parameters.  

\begin{figure}
    \centering
    \begin{subfigure}{0.45\textwidth}
        \centering
        \includegraphics[width=\textwidth]{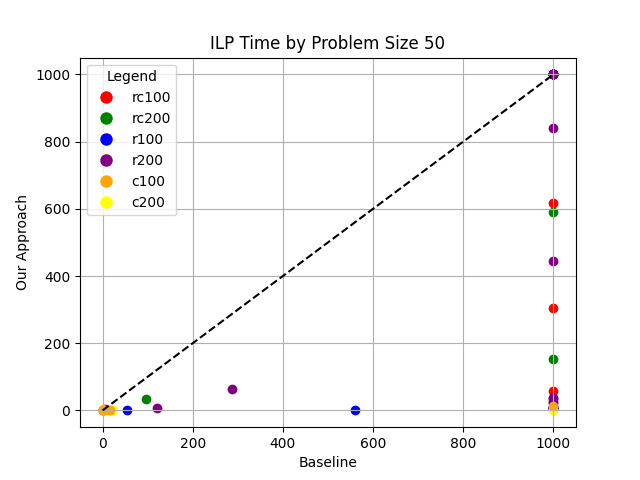}
        \caption{ILP Run Time 50 Customers}
        \label{fig:50runILP}
    \end{subfigure}
    \hfill
    \begin{subfigure}{0.45\textwidth}
        \centering
        \includegraphics[width=\textwidth]{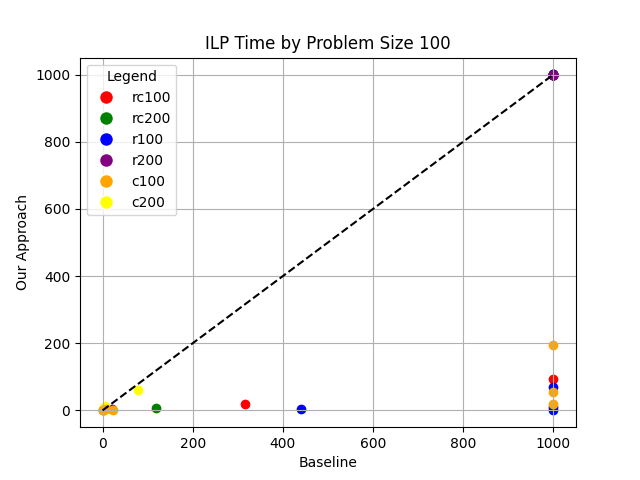}
        \caption{ILP Run Time 100 Customers}
        \label{fig:100runILP}
    \end{subfigure}
    
    \begin{subfigure}{0.45\textwidth}
        \centering
        \includegraphics[width=\textwidth]{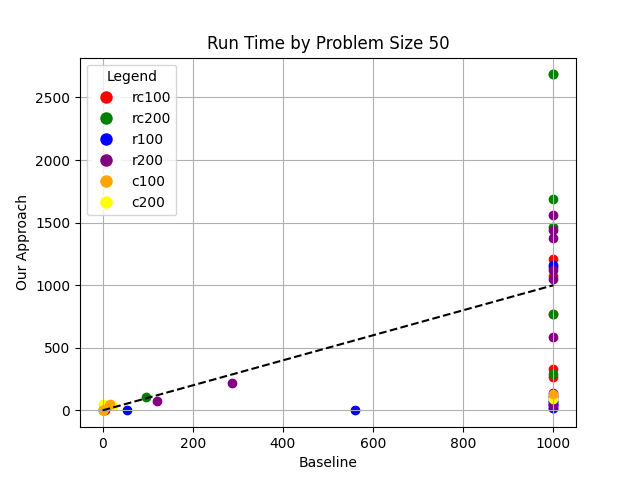}
        \caption{ILP and All LP Total Run Time 50 Customers}
        \label{fig:50run}
    \end{subfigure}
    \hfill
    \begin{subfigure}{0.45\textwidth}
        \centering
        \includegraphics[width=\textwidth]{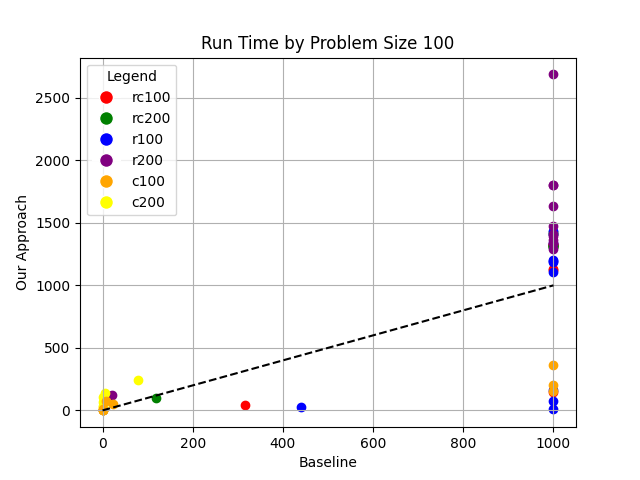}
        \caption{ILP and All LP Total Run Time 100 Customers}
        \label{fig:100run}
    \end{subfigure}
    
    \begin{subfigure}{0.45\textwidth}
        \centering
        \includegraphics[width=\textwidth]{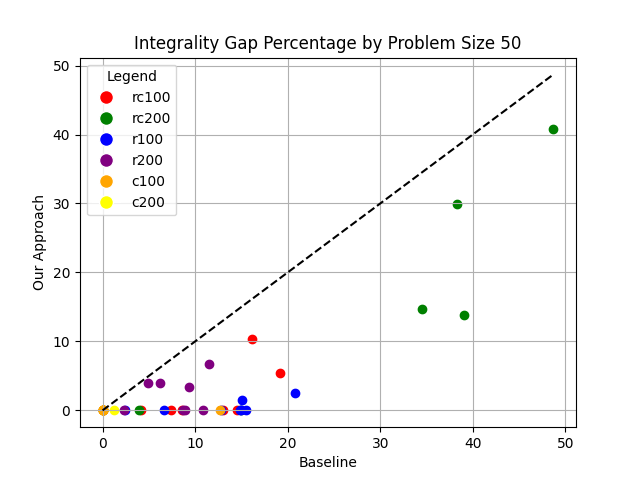}
        \caption{Integrality Gap Percent; 50 Customers}
        \label{fig:50Integ_per}
    \end{subfigure}
    \hfill
    \begin{subfigure}{0.45\textwidth}
        \centering
        \includegraphics[width=\textwidth]{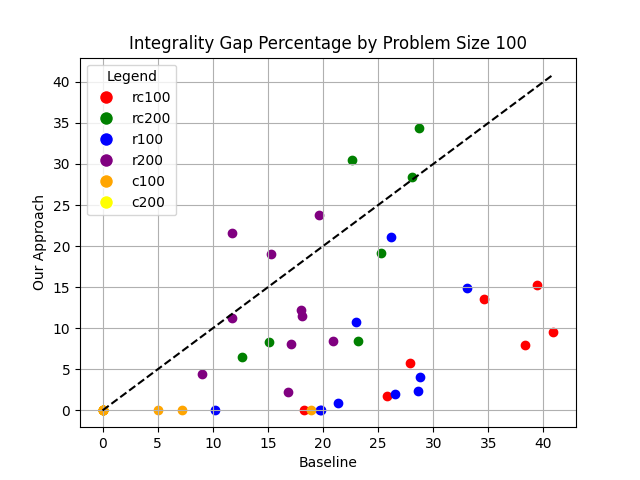}
        \caption{Integrality Gap Percent; 100 Customers}
        \label{fig:100Integ_per}
    \end{subfigure}
    
    \caption{Combined Figure Set One}
    \label{fig:combined}
\end{figure}

\begin{figure}
    \centering
    \begin{subfigure}{0.45\textwidth}
        \centering
        \includegraphics[width=\textwidth]{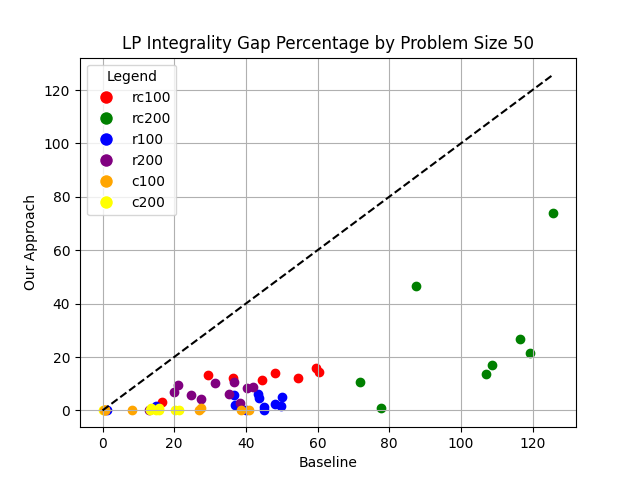}
        \caption{LP Integrality Gap Percent; 50 Customers}
        \label{fig:50LPInteg_per}
    \end{subfigure}
    \hfill
    \begin{subfigure}{0.45\textwidth}
        \centering
        \includegraphics[width=\textwidth]{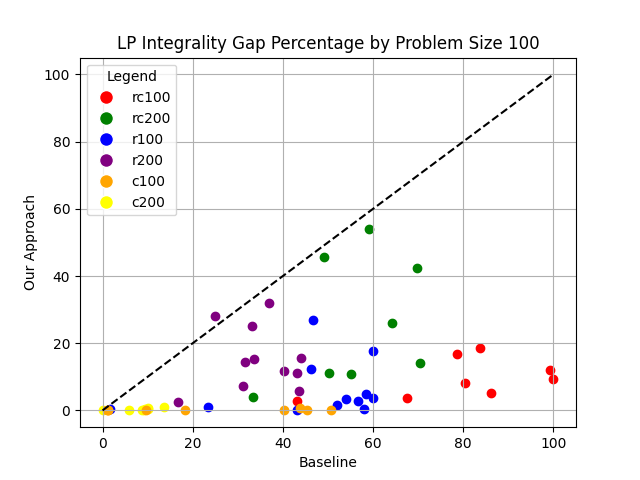}
        \caption{LP Integrality Gap Percent; 100 Customers}
        \label{fig:100LPInteg_per}
    \end{subfigure}
    
    \begin{subfigure}{0.45\textwidth}
        \centering
        \includegraphics[width=\textwidth]{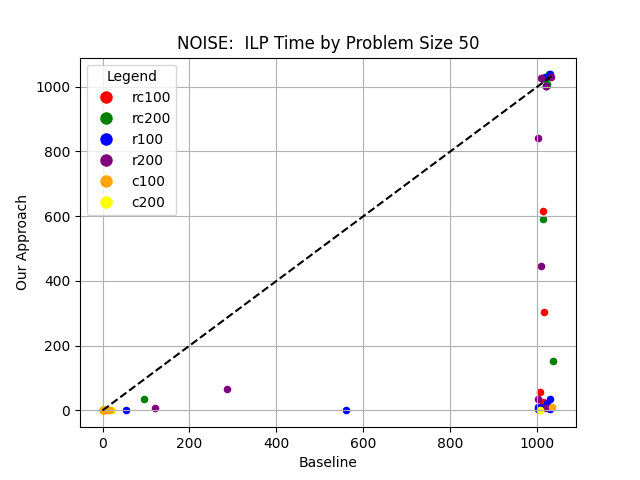}
        \caption{NOISE: ILP Run Time 50 Customers}
        \label{fig:50runILPNOISE}
    \end{subfigure}
    \hfill
    \begin{subfigure}{0.45\textwidth}
        \centering
        \includegraphics[width=\textwidth]{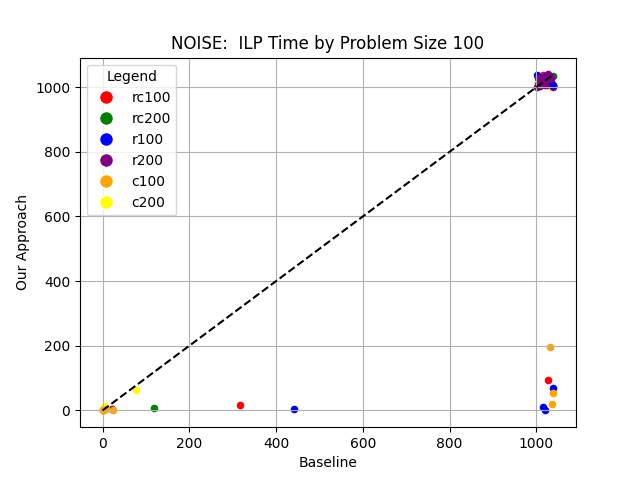}
        \caption{NOISE: ILP Run Time 100 Customers}
        \label{fig:100runILPNOISE}
    \end{subfigure}
    
    \begin{subfigure}{0.45\textwidth}
        \centering
        \includegraphics[width=\textwidth]{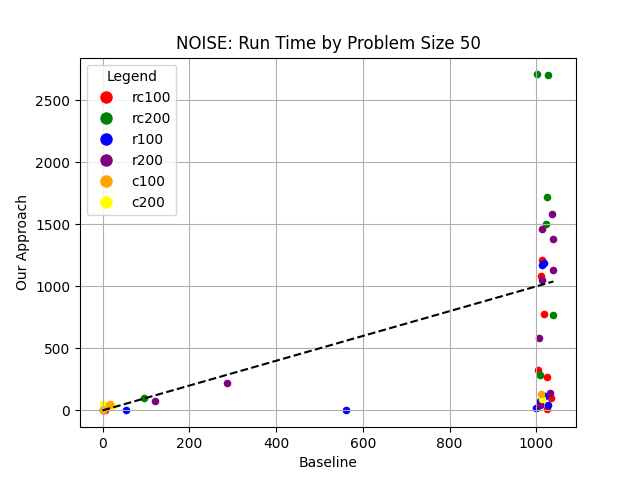}
        \caption{NOISE: ILP and All LP Total Run Time 50 Customers}
        \label{fig:50runNOISE}
    \end{subfigure}
    \hfill
    \begin{subfigure}{0.45\textwidth}
        \centering
        \includegraphics[width=\textwidth]{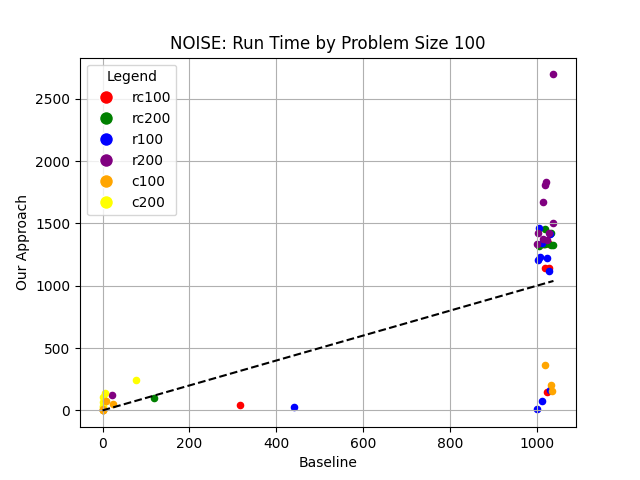}
        \caption{NOISE: ILP and All LP Total Run Time 100 Customers}
        \label{fig:100runNOISE}
    \end{subfigure}
    
    \caption{Combined Figure Set Two}
    \label{fig:combined2}
\end{figure}

We provide  comparisons across Solomon data sets/problem sizes in Figures \ref{fig:50runILP}-\ref{fig:100runILP} (MILP time) and Figures \ref{fig:50run}-\ref{fig:100run} (MILP+ all LPs time) 
Figures \ref{fig:50Integ_per}-\ref{fig:100Integ_per} (integrality gap as a percentage). 
We define the integrality gap to be the difference between the MILP upper bound and the MILP Lower bound which are both produced over the course of running the MILP solver.
For Figures \ref{fig:50Integ_per}-\ref{fig:100Integ_per} we normalize the data as follows.  We take integrality gap and divided it by the MILP upper bound then multiply by 100 to get the value to plot for both baseline and our approach.  
For problems with larger integrality gaps we believe that using valid inequalities to tighten the bound.
In Figures \ref{fig:50LPInteg_per}-\ref{fig:100LPInteg_per} we plot the LP integrality gap.  We define this as the difference between the MILP upper bound and the LP at the root node of the branch \& bound tree.   We normalize the data as follows.  We take LP integrality gap and divided it by the MILP upper bound then multiply by 100 to get the value to plot for both baseline and our approach.  
Since a large number of points for the baseline hit the 1000 second time horizon we had issues with visualization.  Hence for such instances we add random noise to the baseline value  for ILP run time and Total run time in Figs \ref{fig:50runILPNOISE}-\ref{fig:100runNOISE}.

We observe that LA-Discretization performs best relative to the baseline on problem instances which are hardest for baseline MILP to solve efficiently.
%

\pagebreak

\section{Conclusion}
\label{sec_conc}
We have introduced a new approach for the Vehicle Routing Problem with Time Windows (VRPTW) called Local Area Discretization (LA-Discretization). This approach is designed to gain the advantages of the expanded (column generation) formulation for VRPTW \citep{Desrochers1992} and the advantages of the compact two-index formulation without the associated disadvantages. Specifically we gain the advantages of the dramatically tighter LP relaxation of the expanded formulation while being able to solve optimization with an off-the-shelf MILP solver. This advantage should not be underestimated. Unlike related CG approaches no elaborate dual stabilization is required, nor are elaborate pricing algorithms for elementary resource constrained shortest path problems \citep{irnich2005shortest}. However, we find that this method does not always lead to rapid identification of the optimal solution.  In those cases we are close to optimal while the baseline is still far away. 

Our solution is inspired by \citep{boland2017continuous}'s work on time window discretization and pre-computation of component routes localized in space in \citep{mandal2022local_2, mandal2023graph}. Exploiting this approach applies time/capacity window discretization in the LP relaxation, and enforces that the solution to the compact LP is consistent with elementarity locally (hence eliminating sub-tours). LA-Discretization increases and decreases the discretization of time/capacity and the degree of enforcement of elementarity in such a manner as to produce an LP that is compact and tighter than a baseline compact formulation. However the vertices of the associated MILP can express every optimal solution to the original VRPTW. Hence solving this MILP exactly ensures an optimal solution to the original VRPTW.  Utilizing LA-arcs for cycle removal in compact formulations stands out as a significant advancement, offering users the flexibility to scale up their formulations alongside the evolving computational prowess of MILP solvers.

In future work we intend to explore the use of valid inequalities such as subset-row inequalities \citep{jepsen2008subset,wang2017tracking} and rounded capacity inequalities \citep{archetti2011column} for efficient tightening of the LP bound in a manner so as to exploit LA-arcs.  
We also intend to explore the use of  alternative compact LP relaxations with our approach such as that of \citep{gavish1981scheduling}. Furthermore we intend to explore the use of multiple LA-neighborhoods for customers in the same LP.  We also intend to create a graph where nodes correspond to customers and NG  neighbors \citep{baldacci2011new} visited which keeps track of the route in an NG-neighbor like manner so as to further restrict cycles. 

\bibliographystyle{abbrvnat} 
\bibliography{col_gen_bib}

\appendix

\section{Appendix/Extensions Overview}

In this appendix we consider relevant proofs, derivations and extensions to the main document. We organize this appendix as follows. 
In Appendix \ref{sec_proof} we provide a proof certifying the correctness of our approach to generating L-Arcs.  
In Appendix \ref{proofCapBucketFeas} we provide proof that bucket feasibility ensures a sufficient parameterization. In Appendix \ref{sec_proof_converge} we prove that Algorithm \ref{final_opt} converges to a sufficient parameterization. Appendix \ref{allResults} provides additional experimental results directly related to this paper. 


\section{Proof of Equivalent Representation of Departure Time}
\label{sec_proof}

In this section we prove that \eqref{tauDef} accurately characterizes $T_{r}(t)$ for all $r\in R_p,t_0\geq t \geq 0, p \in P^+$.  We prove this using induction. Observe that for the base case where $|N_p|=0$ that the claim holds by definition.  If \eqref{tauDef}  does not hold for all cases then there must exist some $r,t$ s.t \eqref{tauDef} does not hold but \eqref{tauDef} holds for $r^-$ for all time.  We describe this formally below using $u,w$ to denote the first two customers in $r$. 

\textbf{Claim to be Proven False}:  
\begin{align}
\label{claim_to_false}
T_{r}(t)=T_{r^-}(-t_{uw}+\min(t,t^+_u))-\infty[t<t^-_u] \neq -c_r+\min(t,\phi_{r}) -\infty[\min(t,t^+_u)<\hat{\phi}_{r}]
\end{align}
For now we assume \eqref{claim_to_false} is true.  We rew-write \eqref{ituDef} below for reference 
\begin{align}
\label{ituDef2}
    T_{r}(t)=-c_r+\min(t,\phi_{r}) -\infty[\min(t,t^+_{u_p})<\hat{\phi}_{r}] \quad \forall p \in P^+,r \in R^+_p 
\end{align}

\textbf{Proof:  }
We now use \eqref{ituDef} to replace $T_{r^-}(-t_{uw}+\min(t,t^+_u))$ in \eqref{claim_to_false}.  
\begin{align}
\label{eq_term_1}
-c_{r^-}+\min(-t_{uw}+\min(t,t^+_u),\phi_{r^-}) -\infty[-t_{uw}+\min(t,t^+_u)-\infty[t<t^-_u]<\hat{\phi}_{r^-}]-\infty[t<t^-_u]\\
\neq  -c_r+\min(t,\phi_{r}) -\infty[\min(t,t^+_u)<\hat{\phi}_{r}]\nonumber 
\end{align}
We now consider the terms corresponding to $\infty$ and those not corresponding to $\infty$ separately in \eqref{eq_term_1}. For the \eqref{claim_to_false} to be true then it must be the case that  $T_{r^-}(-t_{uw}+\min(t,t^+_u))-\infty[t<t^-_u] \neq -c_r+\min(t,\phi_{r}) -\infty[\min(t,t^+_u)<\hat{\phi}_{r}]$.  A necessary condition for \eqref{claim_to_false} to be true is that  one or both of the following must must be satisfied.  
\begin{subequations}
\label{2Pos}
\begin{align}
 -\infty[-t_{uw}+\min(t,t^+_u)<\hat{\phi}_{r^-}]-\infty[t<t^-_u] \neq -\infty[\min(t,t^+_u)<\hat{\phi}_{r}]\label{pos_1a}\\
 -c_{r^-} +\min(-t_{uw}+\min(t,t^+_u),\phi_{r}) \neq c_r + \min(t,\phi_{r}) \label{pos_2a}
\end{align}
\end{subequations}

Let us consider \eqref{pos_1a} first.  We apply the following transformations including plugging in the definition of $\hat{\phi}_{r}$.  
\begin{subequations}
\begin{align}
\label{infEq}
 -\infty[-t_{uw}+\min(t,t^+_u)<\hat{\phi}_{r^-}]-\infty[t<t^-_u]\neq -\infty[\min(t,t^+_u)<\hat{\phi}_{r}]\\
 -\infty[-t_{uw}+\min(t,t^+_u)<\hat{\phi}_{r^-}]-\infty[t<t^-_u]\neq -\infty[\min(t,t^+_u)<\max(t^-_u,\hat{\phi}_{r^-}+t_{uw})]\\
-\infty[\min(t,t^+_u)<\hat{\phi}_{r^-}+t_{uw}]-\infty[t<t^-_u] \neq -\infty[\min(t,t^+_u)<\max(t^-_u,\hat{\phi}_{r^-}+t_{uw})]\label{finEqAset}
\end{align}
\end{subequations}
Consider that  $t<t^-_u$  then both the LHS and RHS of \eqref{finEqAset} equal $-\infty$.  Thus if a violation occurs in \eqref{finEqAset} then $t\geq t^-_u$.  Thus we gain the following expression for a potential violation:  
\begin{align}
 -\infty[\min(t,t^+_u)<\hat{\phi}_{r^-}+t_{uw}] \neq -\infty[\min(t,t^+_u)<\max(t^-_u,\hat{\phi}_{r^-}+t_{uw})]\label{tryme}
\end{align}
The only way for \eqref{tryme} to be satisfied is if $\hat{\phi}_{r^-}+t_{uw}\leq \min(t,t^+_u)<t^-_u$.  However both $t$ and $t^+_u$ are greater than or equal to $t^-_u$ creating a contradiction.  

We now consider the case where \eqref{pos_2a} is violated.  We now plug in the definition of $\phi_{r}$, into \eqref{pos_2a}, and  noting that $t_{uw}=c_{uw}$ to produce the following expressions by equivalent transformations.  
\begin{subequations}
\begin{align}
-c_{r^-}+\min(-t_{uw}+\min(t,t^+_u),\phi_{r^-})\neq -c_r+\min(t,\phi_{r})\\
 -c_{r^-}+\min(-t_{uw}+\min(t,t^+_u),\phi_{r^-})\neq -c_{r^-}-c_{uw}+\min(t,\min(t^+_u,\phi_{r^-}+t_{uw}))\\
 \min(\min(t-t_{uw},t^+_u-t_{uw}),\phi_{r^-})\neq -c_{uw}+\min(t,\min(t^+_u,\phi_{r^-}+t_{uw}))\\
 \min(\min(t-t_{uw},t^+_u-t_{uw}),\phi_{r^-})\neq -t_{uw}+\min(t,\min(t^+_u,\phi_{r^-}+t_{uw}))\\
\min(t-t_{uw},t^+_u-t_{uw},\phi_{r^-})\neq -t_{uw}+\min(t,t^+_u,\phi_{r^-}+t_{uw})\\
 \min(t-t_{uw},t^+_u-t_{uw},\phi_{r^-}) \neq \min(t-t_{uw},t^+_u-t_{uw},\phi_{r^-})\label{finverGroup}
\end{align}
\end{subequations}

Since the RHS and LHS of \eqref{finverGroup} are identical so the inequality does not hold \eqref{pos_2a} does not hold.  Since neither \eqref{pos_1a} or \eqref{pos_2a} hold then \eqref{claim_to_false} is false.  Thus \eqref{tauDef} must be true.

\section{Bucket Feasibility}
\label{proofCapBucketFeas}

 Consider the following solution $\{N^{\odot}_u,T^{\odot}_u,D^{\odot}_u\} \forall u \in N)$ for $(x,y,z)$ where  $x,y$ are equal to $x^*,y^*$

For each $i,j$ in $E^{D\odot}$ or $E^{T\odot}$ we use $\hat{c}_{ij}=0$ when $u_i=u_j$ and $\hat{c}_{ij}=c_{u_i u_j}$ when $u_i \neq u_j$.   For each $(i,j)$ we use $E^{D\odot}_{ij}$ to denote the edges in a lowest cost path in $E^{D\odot}_{ij}$ starting at $(u_i,d^+_i,d^+_i)$ and ending at $(u_j,d^+_j,d^+_j)$ where edges have cost defined by defined by $\hat{c}$. Similarly for each $(i,j)$ we use $E^{T\odot}_{ij}$ to denote the edges in a lowest cost path in $E^{T\odot}_{ij}$ starting at $(u_i,t^+_i,t^+_i)$ and ending at $(u_j,t^+_j,t^+_j)$ where edges are of cost defined by $\hat{c}$. Notice that for $u_i\neq u_j$ that $E^{D\odot}_{ij}$ contains exactly one edge $\hat{i},\hat{j}$ for which $u_{\hat{i}}=u_i$ and $u_{\hat{j}}=u_j$ and no other edges between non-identical customers. 
\begin{subequations}    
\begin{align}
[u_i=u][u_j=v]=\sum_{\hat{i},\hat{j} \in E^{\odot D}_{ij}}[u_{\hat{i}}=u][u_{\hat{j}}=v] \label{letP1}\\
[u_i=u][u_j=v]=\sum_{\hat{i},\hat{j} \in E^{\odot T}_{ij}}[u_{\hat{i}}=u][u_{\hat{j}}=v]\label{letP2}
\end{align}
\end{subequations}


We now construct a feasible solution to optimization over  given the LP over 

Now initialize $z\leftarrow 0$.  Now iterate over $z^{D*}_{ij}$ and increase $z^D_{\hat{i}\hat{j}}$ by $z^{D*}_{ij}$ for each $\hat{i},\hat{j}$ in $E^{\odot D}_{ij}$.  Now iterate over $z^{T*}_{ij}$ and increase $z^T_{\hat{i}\hat{j}}$ by $z^{T*}_{ij}$ for each $\hat{i},\hat{j}$ in $E^{\odot T}_{ij}$.  

Observe that $x,y$ are feasible for $\{N^*_u,T^{\odot}_u,D^{\odot}_u\} \forall u \in N)$; and since only $x$ is in the objective that we have not changed the objective.  Thus only violations in \eqref{frlow_grdem},\eqref{frlow_grdem2},\eqref{frlow_grtime_1},\eqref{frlow_grtime_2} are possible.  Observe that via \eqref{letP1} that \eqref{frlow_grdem2} is obeyed and via and \eqref{letP2} that \eqref{frlow_grtime_2} is obeyed. Since the intermediate edges in the paths $E^{\odot D}_{ij}$ can be ignored then flow constraints \eqref{frlow_grdem},\eqref{frlow_grtime_1} are obeyed.
 
\section{Proof of Convergence to a Sufficient Parametization}
\label{sec_proof_converge}

In this section we establish that Alg \ref{final_opt} generates a sufficient parameterization when ITER\_MAX=$\infty$. 

Observe that all expansion/contraction operations never decrease the primal objective.  Thus Alg \ref{final_opt} monotonically increases the dual objective.  Thus we need establish that Alg \ref{final_opt} can not converge to a value other than the maximum.  

Observe that after $\zeta$ iterations of not increasing the objective by MIN\_INC that we have $N_u\leftarrow N_u^{\odot}$ for all $u \in N$.  

Given that $N_u\leftarrow N_u^{\odot}$ for all $u \in N$ observe Alg \ref{final_opt} can continue for a maximum of $\sum_{u \in N}(|D_u^{\odot}|+|T_u^{\odot}|)$ iterations before either the bound is increased or the buckets for time /capacity are at minimal granularity;  thus producing a sufficient parameterization.  

\section{All Results}
\label{allResults}

We provide all timing results in this section for our data sets. We provide results on the Solomon instance data sets \citep{solomon1987algorithms}.
Each row describes the performance of the problem instance for a given approach.  The baseline performance on a given instance is the row immediately below the row for LA-Discretization. We use the following columns to describe our results.
\begin{enumerate}
\item File ID:  We provide the problem instance ID. 
\item Approach:  We indicate whether we used LA-Discretization or the Baseline MILP solver.
\item LP Obj:  This describes the root LP objective of the MILP being solved.
\item MILP LB:  
This is the MILP dual bound at the termination of optimization.
\item MILP Obj: This is the MILP objective at the termination of optimization,
\item MILP time: This is the amount of time required to solve the MILP. Note that we cut off optimization at 1000 seconds.  All times are rounded down to the nearest second.  
\item LP time: This is the time spent solving LPs over the course of LA-Discretization. This is the time required to solve the root LP for the baseline.  All times are rounded down to the nearest second.  
\item 10X$<$: This is indicated as YES if LA-Discretization vastly outperforms the baseline for the given instance. We define vastly outperform as meeting the following two criteria (1) solving the optimization problem exactly and (2) doing so either ten times faster than the baseline or if the baseline does not finish during the allotted optimization time.
\end{enumerate}

\pagebreak
\begin{table}
\centering
\caption{Num Customers = 25 Problem Set = rc100}
\label{tab_all:25_rc100}
\begin{tabular}{cccccccccccccccccccccccccccccccccccccccccccccccccccccccc}
\toprule
file num & Approach & lp obj & mip dual bound & ILP obj & ilp time & total lp time & 10X+  \\
\midrule
rc101 & OUR APPROACH & 402.8 & 461.1 & 461.1 & 0.1 & 0.3 &   \\
rc101 & BASELINE MILP & 336.3 & 461.1 & 461.1 & 0.1 & 0.0 &   \\
rc102 & OUR APPROACH & 350.9 & 351.8 & 351.8 & 0.1 & 0.7 &   \\
rc102 & BASELINE MILP & 302.5 & 351.8 & 351.8 & 0.4 & 0.0 &   \\
rc103 & OUR APPROACH & 319.2 & 332.8 & 332.8 & 0.4 & 3.1 &   \\
rc103 & BASELINE MILP & 280.9 & 332.8 & 332.8 & 9.5 & 0.0 &   \\
rc104 & OUR APPROACH & 302.0 & 306.6 & 306.6 & 0.3 & 2.1 &   \\
rc104 & BASELINE MILP & 278.9 & 306.6 & 306.6 & 1.5 & 0.0 &   \\
rc105 & OUR APPROACH & 408.5 & 411.3 & 411.3 & 0.1 & 0.6 &   \\
rc105 & BASELINE MILP & 308.0 & 411.3 & 411.3 & 4.4 & 0.0 &   \\
rc106 & OUR APPROACH & 334.7 & 345.5 & 345.5 & 0.5 & 4.7 &   \\
rc106 & BASELINE MILP & 299.0 & 345.5 & 345.5 & 0.7 & 0.0 &   \\
rc107 & OUR APPROACH & 297.8 & 298.3 & 298.3 & 0.1 & 2.1 &   \\
rc107 & BASELINE MILP & 272.9 & 298.3 & 298.3 & 0.6 & 0.0 &   \\
rc108 & OUR APPROACH & 293.7 & 294.5 & 294.5 & 0.7 & 16.9 &   \\
rc108 & BASELINE MILP & 272.4 & 294.5 & 294.5 & 0.9 & 0.0 &   \\
\bottomrule
\end{tabular}
\end{table}
\begin{table}
\centering
\caption{Num Customers = 25 Problem Set = rc200}
\label{tab_all:25_rc200}
\begin{tabular}{cccccccccccccccccccccccccccccccccccccccccccccccccccccccc}
\toprule
file num & Approach & lp obj & mip dual bound & ILP obj & ilp time & total lp time & 10X+  \\
\midrule
rc201 & OUR APPROACH & 358.8 & 360.2 & 360.2 & 0.0 & 0.2 &   \\
rc201 & BASELINE MILP & 274.0 & 360.2 & 360.2 & 0.1 & 0.0 &   \\
rc202 & OUR APPROACH & 307.3 & 338.0 & 338.0 & 1.5 & 10.1 &   \\
rc202 & BASELINE MILP & 181.7 & 338.0 & 338.0 & 29.0 & 0.0 &   \\
rc203 & OUR APPROACH & 259.6 & 326.9 & 326.9 & 329.4 & 34.6 & YES \\
rc203 & BASELINE MILP & 162.6 & 245.7 & 326.9 & 1000.0 & 0.0 &    \\
rc204 & OUR APPROACH & 238.0 & 269.0 & 299.7 & 1000.1 & 11.8 &   \\
rc204 & BASELINE MILP & 157.8 & 195.3 & 299.7 & 1000.0 & 0.0 &   \\
rc205 & OUR APPROACH & 314.1 & 338.0 & 338.0 & 0.1 & 1.1 &   \\
rc205 & BASELINE MILP & 193.7 & 338.0 & 338.0 & 0.4 & 0.0 &   \\
rc206 & OUR APPROACH & 301.3 & 324.0 & 324.0 & 0.6 & 2.8 &   \\
rc206 & BASELINE MILP & 183.2 & 324.0 & 324.0 & 1.7 & 0.0 &   \\
rc207 & OUR APPROACH & 254.5 & 298.3 & 298.3 & 34.8 & 14.4 & YES \\
rc207 & BASELINE MILP & 153.7 & 274.2 & 298.3 & 1000.0 & 0.0 &    \\
rc208 & OUR APPROACH & 168.1 & 186.0 & 269.1 & 1000.0 & 17.4 &   \\
rc208 & BASELINE MILP & 148.3 & 190.0 & 269.1 & 1000.0 & 0.0 &   \\
\bottomrule
\end{tabular}
\end{table}
\begin{table}
\centering
\caption{Num Customers = 25 Problem Set = r100}
\label{tab_all:25_r100}
\begin{tabular}{cccccccccccccccccccccccccccccccccccccccccccccccccccccccc}
\toprule
file num & Approach & lp obj & mip dual bound & ILP obj & ilp time & total lp time & 10X+  \\
\midrule
r101 & OUR APPROACH & 617.1 & 617.1 & 617.1 & 0.0 & 0.0 &   \\
r101 & BASELINE MILP & 617.1 & 617.1 & 617.1 & 0.0 & 0.0 &   \\
r102 & OUR APPROACH & 546.3 & 547.1 & 547.1 & 0.0 & 0.1 &   \\
r102 & BASELINE MILP & 384.7 & 547.1 & 547.1 & 0.4 & 0.0 &   \\
r103 & OUR APPROACH & 454.2 & 454.6 & 454.6 & 0.1 & 0.9 &   \\
r103 & BASELINE MILP & 328.9 & 454.6 & 454.6 & 3.9 & 0.0 &   \\
r104 & OUR APPROACH & 415.0 & 416.9 & 416.9 & 0.1 & 1.4 & YES \\
r104 & BASELINE MILP & 319.5 & 416.9 & 416.9 & 36.2 & 0.0 &    \\
r105 & OUR APPROACH & 530.5 & 530.5 & 530.5 & 0.0 & 0.0 &   \\
r105 & BASELINE MILP & 456.5 & 530.5 & 530.5 & 0.1 & 0.0 &   \\
r106 & OUR APPROACH & 457.3 & 465.4 & 465.4 & 0.1 & 0.3 &   \\
r106 & BASELINE MILP & 346.3 & 465.4 & 465.4 & 0.6 & 0.0 &   \\
r107 & OUR APPROACH & 417.4 & 424.3 & 424.3 & 0.1 & 1.8 &   \\
r107 & BASELINE MILP & 319.7 & 424.3 & 424.3 & 6.8 & 0.0 &   \\
r108 & OUR APPROACH & 391.7 & 397.3 & 397.3 & 0.2 & 5.5 &   \\
r108 & BASELINE MILP & 311.0 & 397.3 & 397.3 & 20.4 & 0.0 &   \\
r109 & OUR APPROACH & 440.3 & 441.3 & 441.3 & 0.1 & 0.4 &   \\
r109 & BASELINE MILP & 347.3 & 441.3 & 441.3 & 0.2 & 0.0 &   \\
r110 & OUR APPROACH & 421.5 & 444.1 & 444.1 & 1.0 & 2.7 & YES \\
r110 & BASELINE MILP & 312.8 & 444.1 & 444.1 & 65.8 & 0.0 &    \\
r111 & OUR APPROACH & 418.4 & 428.8 & 428.8 & 0.2 & 1.8 &   \\
r111 & BASELINE MILP & 320.6 & 428.8 & 428.8 & 5.3 & 0.0 &   \\
r112 & OUR APPROACH & 371.3 & 393.0 & 393.0 & 2.2 & 9.7 & YES \\
r112 & BASELINE MILP & 309.9 & 393.0 & 393.0 & 353.0 & 0.0 &    \\
\bottomrule
\end{tabular}
\end{table}
\begin{table}
\centering
\caption{Num Customers = 25 Problem Set = r200}
\label{tab_all:25_r200}
\begin{tabular}{cccccccccccccccccccccccccccccccccccccccccccccccccccccccc}
\toprule
file num & Approach & lp obj & mip dual bound & ILP obj & ilp time & total lp time & 10X+  \\
\midrule
r201 & OUR APPROACH & 458.6 & 463.3 & 463.3 & 0.0 & 0.4 &   \\
r201 & BASELINE MILP & 417.0 & 463.3 & 463.3 & 0.0 & 0.0 &   \\
r202 & OUR APPROACH & 399.6 & 410.5 & 410.5 & 0.2 & 2.0 &   \\
r202 & BASELINE MILP & 314.7 & 410.5 & 410.5 & 1.5 & 0.0 &   \\
r203 & OUR APPROACH & 362.5 & 391.4 & 391.4 & 0.9 & 8.8 &   \\
r203 & BASELINE MILP & 302.3 & 391.4 & 391.4 & 43.0 & 0.0 &   \\
r204 & OUR APPROACH & 322.2 & 355.0 & 355.0 & 2.5 & 15.7 &   \\
r204 & BASELINE MILP & 292.6 & 355.0 & 355.0 & 29.6 & 0.0 &   \\
r205 & OUR APPROACH & 386.4 & 393.0 & 393.0 & 0.1 & 2.4 &   \\
r205 & BASELINE MILP & 337.9 & 393.0 & 393.0 & 0.6 & 0.0 &   \\
r206 & OUR APPROACH & 349.2 & 374.4 & 374.4 & 1.7 & 11.3 &   \\
r206 & BASELINE MILP & 292.9 & 374.4 & 374.4 & 5.1 & 0.0 &   \\
r207 & OUR APPROACH & 338.4 & 361.6 & 361.6 & 2.9 & 23.3 &   \\
r207 & BASELINE MILP & 292.4 & 361.6 & 361.6 & 11.2 & 0.0 &   \\
r208 & OUR APPROACH & 312.3 & 328.2 & 328.2 & 2.0 & 8.1 &   \\
r208 & BASELINE MILP & 292.3 & 328.2 & 328.2 & 2.9 & 0.0 &   \\
r209 & OUR APPROACH & 354.6 & 370.7 & 370.7 & 0.9 & 5.3 &   \\
r209 & BASELINE MILP & 303.3 & 370.7 & 370.7 & 0.8 & 0.0 &   \\
r210 & OUR APPROACH & 378.3 & 404.6 & 404.6 & 1.9 & 9.5 &   \\
r210 & BASELINE MILP & 302.7 & 404.6 & 404.6 & 1.8 & 0.0 &   \\
r211 & OUR APPROACH & 313.6 & 350.9 & 350.9 & 19.2 & 19.1 &   \\
r211 & BASELINE MILP & 292.4 & 350.9 & 350.9 & 44.7 & 0.0 &   \\
\bottomrule
\end{tabular}
\end{table}
\begin{table}
\centering
\caption{Num Customers = 25 Problem Set = c100}
\label{tab_all:25_c100}
\begin{tabular}{cccccccccccccccccccccccccccccccccccccccccccccccccccccccc}
\toprule
file num & Approach & lp obj & mip dual bound & ILP obj & ilp time & total lp time & 10X+  \\
\midrule
c101 & OUR APPROACH & 191.3 & 191.3 & 191.3 & 0.0 & 0.1 &   \\
c101 & BASELINE MILP & 191.3 & 191.3 & 191.3 & 0.0 & 0.0 &   \\
c102 & OUR APPROACH & 190.3 & 190.3 & 190.3 & 0.0 & 0.4 &   \\
c102 & BASELINE MILP & 172.7 & 190.3 & 190.3 & 0.2 & 0.0 &   \\
c103 & OUR APPROACH & 190.3 & 190.3 & 190.3 & 0.1 & 2.0 &   \\
c103 & BASELINE MILP & 163.8 & 190.3 & 190.3 & 0.6 & 0.0 &   \\
c104 & OUR APPROACH & 186.9 & 186.9 & 186.9 & 0.1 & 2.8 &   \\
c104 & BASELINE MILP & 160.5 & 186.9 & 186.9 & 2.4 & 0.0 &   \\
c105 & OUR APPROACH & 191.3 & 191.3 & 191.3 & 0.0 & 0.1 &   \\
c105 & BASELINE MILP & 191.2 & 191.3 & 191.3 & 0.0 & 0.0 &   \\
c106 & OUR APPROACH & 191.3 & 191.3 & 191.3 & 0.0 & 0.1 &   \\
c106 & BASELINE MILP & 191.3 & 191.3 & 191.3 & 0.0 & 0.0 &   \\
c107 & OUR APPROACH & 191.3 & 191.3 & 191.3 & 0.0 & 0.1 &   \\
c107 & BASELINE MILP & 191.2 & 191.3 & 191.3 & 0.0 & 0.0 &   \\
c108 & OUR APPROACH & 191.3 & 191.3 & 191.3 & 0.0 & 1.1 &   \\
c108 & BASELINE MILP & 131.5 & 191.3 & 191.3 & 0.2 & 0.0 &   \\
c109 & OUR APPROACH & 190.6 & 191.3 & 191.3 & 0.2 & 6.4 &   \\
c109 & BASELINE MILP & 131.1 & 191.3 & 191.3 & 0.8 & 0.0 &   \\
\bottomrule
\end{tabular}
\end{table}
\begin{table}
\centering
\caption{Num Customers = 25 Problem Set = c200}
\label{tab_all:25_c200}
\begin{tabular}{cccccccccccccccccccccccccccccccccccccccccccccccccccccccc}
\toprule
file num & Approach & lp obj & mip dual bound & ILP obj & ilp time & total lp time & 10X+  \\
\midrule
c201 & OUR APPROACH & 214.7 & 214.7 & 214.7 & 0.0 & 0.1 &   \\
c201 & BASELINE MILP & 214.7 & 214.7 & 214.7 & 0.0 & 0.0 &   \\
c202 & OUR APPROACH & 214.7 & 214.7 & 214.7 & 0.0 & 2.0 &   \\
c202 & BASELINE MILP & 175.0 & 214.7 & 214.7 & 0.6 & 0.0 &   \\
c203 & OUR APPROACH & 214.7 & 214.7 & 214.7 & 0.0 & 3.3 &   \\
c203 & BASELINE MILP & 169.5 & 214.7 & 214.7 & 1.5 & 0.0 &   \\
c204 & OUR APPROACH & 209.6 & 213.1 & 213.1 & 1.0 & 8.6 &   \\
c204 & BASELINE MILP & 165.4 & 213.1 & 213.1 & 24.9 & 0.0 &   \\
c205 & OUR APPROACH & 214.7 & 214.7 & 214.7 & 0.0 & 0.3 &   \\
c205 & BASELINE MILP & 176.4 & 214.7 & 214.7 & 0.0 & 0.0 &   \\
c206 & OUR APPROACH & 214.7 & 214.7 & 214.7 & 0.0 & 0.5 &   \\
c206 & BASELINE MILP & 173.1 & 214.7 & 214.7 & 0.2 & 0.0 &   \\
c207 & OUR APPROACH & 211.3 & 214.5 & 214.5 & 0.2 & 5.8 &   \\
c207 & BASELINE MILP & 166.7 & 214.5 & 214.5 & 1.0 & 0.0 &   \\
c208 & OUR APPROACH & 211.7 & 214.5 & 214.5 & 0.1 & 1.1 &   \\
c208 & BASELINE MILP & 166.4 & 214.5 & 214.5 & 0.3 & 0.0 &   \\
\bottomrule
\end{tabular}
\end{table}
\begin{table}
\centering
\caption{Num Customers = 50 Problem Set = rc100}
\label{tab_all:50_rc100}
\begin{tabular}{cccccccccccccccccccccccccccccccccccccccccccccccccccccccc}
\toprule
file num & Approach & lp obj & mip dual bound & ILP obj & ilp time & total lp time & 10X+  \\
\midrule
rc101 & OUR APPROACH & 841.4 & 944.0 & 944.0 & 2.7 & 3.8 &   \\
rc101 & BASELINE MILP & 610.6 & 944.0 & 944.0 & 5.3 & 0.0 &   \\
rc102 & OUR APPROACH & 710.2 & 822.5 & 822.5 & 303.9 & 24.7 & YES \\
rc102 & BASELINE MILP & 515.9 & 751.5 & 822.5 & 1000.0 & 0.0 &    \\
rc103 & OUR APPROACH & 623.0 & 672.2 & 710.9 & 1000.0 & 73.6 &   \\
rc103 & BASELINE MILP & 480.0 & 575.3 & 711.5 & 1000.0 & 0.0 &   \\
rc104 & OUR APPROACH & 529.8 & 545.8 & 545.8 & 56.9 & 208.9 & YES \\
rc104 & BASELINE MILP & 467.6 & 523.2 & 545.8 & 1000.0 & 0.0 &    \\
rc105 & OUR APPROACH & 747.9 & 855.3 & 855.3 & 6.8 & 8.5 & YES \\
rc105 & BASELINE MILP & 533.6 & 791.8 & 855.3 & 1000.1 & 0.0 &    \\
rc106 & OUR APPROACH & 648.5 & 723.2 & 723.2 & 27.1 & 69.7 & YES \\
rc106 & BASELINE MILP & 500.9 & 618.3 & 723.2 & 999.9 & 0.0 &    \\
rc107 & OUR APPROACH & 573.3 & 642.7 & 642.7 & 616.1 & 157.3 & YES \\
rc107 & BASELINE MILP & 470.9 & 559.0 & 642.7 & 1000.0 & 0.0 &    \\
rc108 & OUR APPROACH & 528.4 & 536.2 & 598.1 & 1000.0 & 208.0 &   \\
rc108 & BASELINE MILP & 463.2 & 502.7 & 599.4 & 1000.0 & 0.0 &   \\
\bottomrule
\end{tabular}
\end{table}
\begin{table}
\centering
\caption{Num Customers = 50 Problem Set = rc200}
\label{tab_all:50_rc200}
\begin{tabular}{cccccccccccccccccccccccccccccccccccccccccccccccccccccccc}
\toprule
file num & Approach & lp obj & mip dual bound & ILP obj & ilp time & total lp time & 10X+  \\
\midrule
rc201 & OUR APPROACH & 677.9 & 684.8 & 684.8 & 0.2 & 2.0 &   \\
rc201 & BASELINE MILP & 385.5 & 684.8 & 684.8 & 0.8 & 0.0 &   \\
rc202 & OUR APPROACH & 539.8 & 613.6 & 613.6 & 591.1 & 180.3 & YES \\
rc202 & BASELINE MILP & 296.3 & 522.1 & 613.6 & 1000.0 & 0.0 &    \\
rc203 & OUR APPROACH & 437.6 & 474.0 & 555.3 & 1000.0 & 463.5 &   \\
rc203 & BASELINE MILP & 257.3 & 365.2 & 557.3 & 1000.0 & 0.0 &   \\
rc204 & OUR APPROACH & 307.0 & 315.1 & 449.8 & 1000.0 & 1683.8 &   \\
rc204 & BASELINE MILP & 237.0 & 273.9 & 444.2 & 1000.0 & 0.0 &   \\
rc205 & OUR APPROACH & 569.4 & 630.2 & 630.2 & 34.1 & 69.3 &   \\
rc205 & BASELINE MILP & 366.8 & 630.2 & 630.2 & 95.5 & 0.0 &   \\
rc206 & OUR APPROACH & 521.2 & 610.0 & 610.0 & 153.6 & 133.2 & YES \\
rc206 & BASELINE MILP & 292.2 & 586.0 & 610.0 & 1000.0 & 0.0 &    \\
rc207 & OUR APPROACH & 462.5 & 483.5 & 561.5 & 1000.0 & 692.5 &   \\
rc207 & BASELINE MILP & 255.3 & 341.7 & 560.2 & 1000.0 & 0.0 &   \\
rc208 & OUR APPROACH & 283.8 & 292.5 & 493.8 & 1000.2 & 1682.9 &   \\
rc208 & BASELINE MILP & 229.2 & 265.5 & 517.7 & 1000.0 & 0.0 &   \\
\bottomrule
\end{tabular}
\end{table}
\begin{table}
\centering
\caption{Num Customers = 50 Problem Set = r100}
\label{tab_all:50_r100}
\begin{tabular}{cccccccccccccccccccccccccccccccccccccccccccccccccccccccc}
\toprule
file num & Approach & lp obj & mip dual bound & ILP obj & ilp time & total lp time & 10X+  \\
\midrule
r101 & OUR APPROACH & 1044.0 & 1044.0 & 1044.0 & 0.0 & 0.1 &   \\
r101 & BASELINE MILP & 1029.8 & 1044.0 & 1044.0 & 0.0 & 0.0 &   \\
r102 & OUR APPROACH & 909.0 & 909.0 & 909.0 & 0.1 & 1.2 & YES \\
r102 & BASELINE MILP & 649.7 & 909.0 & 909.0 & 54.9 & 0.0 &    \\
r103 & OUR APPROACH & 762.1 & 772.9 & 772.9 & 3.1 & 12.7 & YES \\
r103 & BASELINE MILP & 533.1 & 721.7 & 772.9 & 1000.0 & 0.0 &    \\
r104 & OUR APPROACH & 612.1 & 625.4 & 625.4 & 22.7 & 95.3 & YES \\
r104 & BASELINE MILP & 464.6 & 541.6 & 636.0 & 1000.0 & 0.0 &    \\
r105 & OUR APPROACH & 898.5 & 911.8 & 911.8 & 0.3 & 0.8 &   \\
r105 & BASELINE MILP & 793.0 & 911.8 & 911.8 & 1.2 & 0.0 &   \\
r106 & OUR APPROACH & 792.7 & 793.0 & 793.0 & 0.3 & 4.1 & YES \\
r106 & BASELINE MILP & 547.0 & 793.0 & 793.0 & 561.0 & 0.0 &    \\
r107 & OUR APPROACH & 700.5 & 711.1 & 711.1 & 4.4 & 37.1 & YES \\
r107 & BASELINE MILP & 483.2 & 631.7 & 724.4 & 1000.0 & 0.0 &    \\
r108 & OUR APPROACH & 584.9 & 608.1 & 617.7 & 1000.0 & 160.2 &   \\
r108 & BASELINE MILP & 457.5 & 530.6 & 624.7 & 1000.0 & 0.0 &   \\
r109 & OUR APPROACH & 752.1 & 786.8 & 786.8 & 9.5 & 37.9 & YES \\
r109 & BASELINE MILP & 547.2 & 767.5 & 786.8 & 1000.0 & 0.0 &    \\
r110 & OUR APPROACH & 681.1 & 697.0 & 697.0 & 9.9 & 68.5 & YES \\
r110 & BASELINE MILP & 481.0 & 606.3 & 713.0 & 1000.0 & 0.0 &    \\
r111 & OUR APPROACH & 673.0 & 707.2 & 707.2 & 35.5 & 31.5 & YES \\
r111 & BASELINE MILP & 475.4 & 603.8 & 714.2 & 1000.0 & 0.0 &    \\
r112 & OUR APPROACH & 593.3 & 613.9 & 630.2 & 1000.0 & 146.8 &   \\
r112 & BASELINE MILP & 451.8 & 513.2 & 647.9 & 1000.0 & 0.0 &   \\
\bottomrule
\end{tabular}
\end{table}
\begin{table}
\centering
\caption{Num Customers = 50 Problem Set = r200}
\label{tab_all:50_r200}
\begin{tabular}{cccccccccccccccccccccccccccccccccccccccccccccccccccccccc}
\toprule
file num & Approach & lp obj & mip dual bound & ILP obj & ilp time & total lp time & 10X+  \\
\midrule
r201 & OUR APPROACH & 789.5 & 791.9 & 791.9 & 0.1 & 3.5 &   \\
r201 & BASELINE MILP & 700.9 & 791.9 & 791.9 & 0.3 & 0.0 &   \\
r202 & OUR APPROACH & 680.5 & 698.5 & 698.5 & 6.6 & 39.4 & YES \\
r202 & BASELINE MILP & 504.4 & 682.4 & 698.5 & 1000.0 & 0.0 &    \\
r203 & OUR APPROACH & 570.2 & 605.3 & 605.3 & 36.1 & 104.7 & YES \\
r203 & BASELINE MILP & 447.9 & 540.0 & 605.9 & 1000.0 & 0.0 &    \\
r204 & OUR APPROACH & 464.6 & 487.8 & 508.1 & 1000.0 & 378.8 &   \\
r204 & BASELINE MILP & 420.4 & 477.4 & 508.9 & 1000.0 & 0.0 &   \\
r205 & OUR APPROACH & 662.2 & 690.1 & 690.1 & 6.7 & 65.9 &   \\
r205 & BASELINE MILP & 540.9 & 690.1 & 690.1 & 120.7 & 0.0 &   \\
r206 & OUR APPROACH & 584.4 & 632.4 & 632.4 & 445.2 & 141.4 & YES \\
r206 & BASELINE MILP & 455.2 & 581.8 & 639.0 & 1000.0 & 0.0 &    \\
r207 & OUR APPROACH & 519.9 & 556.5 & 576.1 & 1000.0 & 438.9 &   \\
r207 & BASELINE MILP & 420.6 & 521.6 & 575.5 & 1000.0 & 0.0 &   \\
r208 & OUR APPROACH & 455.6 & 468.1 & 487.7 & 1000.0 & 118.7 &   \\
r208 & BASELINE MILP & 412.1 & 469.9 & 494.2 & 1000.0 & 0.0 &   \\
r209 & OUR APPROACH & 568.6 & 600.6 & 600.6 & 64.7 & 153.6 &   \\
r209 & BASELINE MILP & 481.2 & 600.6 & 600.6 & 287.4 & 0.0 &   \\
r210 & OUR APPROACH & 592.8 & 645.6 & 645.6 & 841.9 & 211.2 & YES \\
r210 & BASELINE MILP & 456.2 & 591.2 & 647.8 & 1000.0 & 0.0 &    \\
r211 & OUR APPROACH & 489.8 & 504.6 & 540.9 & 1000.0 & 563.4 &   \\
r211 & BASELINE MILP & 412.1 & 479.3 & 541.3 & 1000.0 & 0.0 &   \\
\bottomrule
\end{tabular}
\end{table}
\begin{table}
\centering
\caption{Num Customers = 50 Problem Set = c100}
\label{tab_all:50_c100}
\begin{tabular}{cccccccccccccccccccccccccccccccccccccccccccccccccccccccc}
\toprule
file num & Approach & lp obj & mip dual bound & ILP obj & ilp time & total lp time & 10X+  \\
\midrule
c101 & OUR APPROACH & 362.4 & 362.4 & 362.4 & 0.1 & 0.5 &   \\
c101 & BASELINE MILP & 361.4 & 362.4 & 362.4 & 0.0 & 0.0 &   \\
c102 & OUR APPROACH & 361.4 & 361.4 & 361.4 & 0.2 & 2.4 &   \\
c102 & BASELINE MILP & 334.0 & 361.4 & 361.4 & 0.5 & 0.0 &   \\
c103 & OUR APPROACH & 361.4 & 361.4 & 361.4 & 0.8 & 40.3 &   \\
c103 & BASELINE MILP & 284.7 & 361.4 & 361.4 & 15.1 & 0.0 &   \\
c104 & OUR APPROACH & 355.2 & 358.0 & 358.0 & 11.9 & 118.5 & YES \\
c104 & BASELINE MILP & 281.0 & 312.4 & 358.0 & 1000.0 & 0.0 &    \\
c105 & OUR APPROACH & 362.4 & 362.4 & 362.4 & 0.1 & 1.2 &   \\
c105 & BASELINE MILP & 360.0 & 362.4 & 362.4 & 0.1 & 0.0 &   \\
c106 & OUR APPROACH & 362.4 & 362.4 & 362.4 & 0.1 & 0.9 &   \\
c106 & BASELINE MILP & 360.1 & 362.4 & 362.4 & 0.1 & 0.0 &   \\
c107 & OUR APPROACH & 362.4 & 362.4 & 362.4 & 0.1 & 1.7 &   \\
c107 & BASELINE MILP & 360.0 & 362.4 & 362.4 & 0.1 & 0.0 &   \\
c108 & OUR APPROACH & 362.4 & 362.4 & 362.4 & 0.2 & 6.2 &   \\
c108 & BASELINE MILP & 261.1 & 362.4 & 362.4 & 0.9 & 0.0 &   \\
c109 & OUR APPROACH & 361.7 & 362.4 & 362.4 & 1.8 & 51.0 &   \\
c109 & BASELINE MILP & 257.1 & 362.4 & 362.4 & 16.1 & 0.0 &   \\
\bottomrule
\end{tabular}
\end{table}
\begin{table}
\centering
\caption{Num Customers = 50 Problem Set = c200}
\label{tab_all:50_c200}
\begin{tabular}{cccccccccccccccccccccccccccccccccccccccccccccccccccccccc}
\toprule
file num & Approach & lp obj & mip dual bound & ILP obj & ilp time & total lp time & 10X+  \\
\midrule
c201 & OUR APPROACH & 360.2 & 360.2 & 360.2 & 0.0 & 0.9 &   \\
c201 & BASELINE MILP & 360.2 & 360.2 & 360.2 & 0.0 & 0.0 &   \\
c202 & OUR APPROACH & 360.2 & 360.2 & 360.2 & 0.1 & 11.9 &   \\
c202 & BASELINE MILP & 318.2 & 360.2 & 360.2 & 2.0 & 0.0 &   \\
c203 & OUR APPROACH & 359.8 & 359.8 & 359.8 & 0.2 & 35.5 &   \\
c203 & BASELINE MILP & 299.1 & 359.8 & 359.8 & 22.5 & 0.0 &   \\
c204 & OUR APPROACH & 349.9 & 350.1 & 350.1 & 0.8 & 89.8 & YES \\
c204 & BASELINE MILP & 288.6 & 345.8 & 350.1 & 1000.0 & 0.0 &    \\
c205 & OUR APPROACH & 359.8 & 359.8 & 359.8 & 0.1 & 2.5 &   \\
c205 & BASELINE MILP & 313.1 & 359.8 & 359.8 & 0.2 & 0.0 &   \\
c206 & OUR APPROACH & 359.8 & 359.8 & 359.8 & 0.1 & 4.4 &   \\
c206 & BASELINE MILP & 310.9 & 359.8 & 359.8 & 0.9 & 0.0 &   \\
c207 & OUR APPROACH & 356.5 & 359.6 & 359.6 & 1.8 & 49.4 &   \\
c207 & BASELINE MILP & 310.7 & 359.6 & 359.6 & 1.3 & 0.0 &   \\
c208 & OUR APPROACH & 347.4 & 350.5 & 350.5 & 0.2 & 7.7 &   \\
c208 & BASELINE MILP & 309.0 & 350.5 & 350.5 & 0.8 & 0.0 &   \\
\bottomrule
\end{tabular}
\end{table}
\begin{table}
\centering
\caption{Num Customers = 100 Problem Set = rc100}
\label{tab_all:100_rc100}
\begin{tabular}{cccccccccccccccccccccccccccccccccccccccccccccccccccccccc}
\toprule
file num & Approach & lp obj & mip dual bound & ILP obj & ilp time & total lp time & 10X+  \\
\midrule
rc101 & OUR APPROACH & 1577.3 & 1619.8 & 1619.8 & 17.8 & 28.8 &   \\
rc101 & BASELINE MILP & 1131.3 & 1619.8 & 1619.8 & 316.1 & 0.1 &   \\
rc102 & OUR APPROACH & 1384.2 & 1431.9 & 1457.4 & 1000.1 & 126.4 &   \\
rc102 & BASELINE MILP & 795.6 & 1098.7 & 1481.6 & 1000.0 & 0.0 &   \\
rc103 & OUR APPROACH & 1176.8 & 1185.9 & 1288.1 & 1000.0 & 320.6 &   \\
rc103 & BASELINE MILP & 694.5 & 857.3 & 1389.1 & 1000.1 & 0.1 &   \\
rc104 & OUR APPROACH & 1050.6 & 1060.8 & 1227.4 & 1000.0 & 327.8 &   \\
rc104 & BASELINE MILP & 670.4 & 782.6 & 1197.0 & 1000.1 & 0.1 &   \\
rc105 & OUR APPROACH & 1461.4 & 1513.7 & 1513.7 & 94.9 & 52.3 & YES \\
rc105 & BASELINE MILP & 927.5 & 1268.4 & 1552.8 & 1000.0 & 0.0 &    \\
rc106 & OUR APPROACH & 1269.2 & 1294.2 & 1373.5 & 1000.1 & 133.0 &   \\
rc106 & BASELINE MILP & 782.3 & 1018.0 & 1411.6 & 1000.1 & 0.0 &   \\
rc107 & OUR APPROACH & 1125.0 & 1139.8 & 1260.9 & 1000.1 & 333.8 &   \\
rc107 & BASELINE MILP & 681.9 & 803.0 & 1358.8 & 1000.1 & 0.1 &   \\
rc108 & OUR APPROACH & 1041.3 & 1047.3 & 1235.9 & 1000.1 & 319.9 &   \\
rc108 & BASELINE MILP & 663.1 & 738.0 & 1217.6 & 999.4 & 0.1 &   \\
\bottomrule
\end{tabular}
\end{table}
\begin{table}
\centering
\caption{Num Customers = 100 Problem Set = rc200}
\label{tab_all:100_rc200}
\begin{tabular}{cccccccccccccccccccccccccccccccccccccccccccccccccccccccc}
\toprule
file num & Approach & lp obj & mip dual bound & ILP obj & ilp time & total lp time & 10X+  \\
\midrule
rc201 & OUR APPROACH & 1214.0 & 1261.8 & 1261.8 & 6.8 & 90.1 &   \\
rc201 & BASELINE MILP & 946.2 & 1261.8 & 1261.8 & 119.0 & 0.0 &   \\
rc202 & OUR APPROACH & 963.6 & 1008.1 & 1100.7 & 1000.4 & 428.7 &   \\
rc202 & BASELINE MILP & 659.9 & 863.3 & 1124.3 & 1000.0 & 0.0 &   \\
rc203 & OUR APPROACH & 764.0 & 779.5 & 1088.5 & 1000.0 & 305.4 &   \\
rc203 & BASELINE MILP & 575.8 & 702.1 & 976.7 & 1000.1 & 0.1 &   \\
rc204 & OUR APPROACH & 647.3 & 655.4 & 943.0 & 1000.1 & 312.3 &   \\
rc204 & BASELINE MILP & 544.5 & 628.9 & 812.4 & 1000.0 & 0.1 &   \\
rc205 & OUR APPROACH & 1045.5 & 1083.8 & 1159.2 & 1000.0 & 317.5 &   \\
rc205 & BASELINE MILP & 759.1 & 1029.0 & 1178.1 & 1000.0 & 0.0 &   \\
rc206 & OUR APPROACH & 951.1 & 969.7 & 1058.1 & 1000.0 & 414.7 &   \\
rc206 & BASELINE MILP & 710.3 & 905.4 & 1066.7 & 1000.0 & 0.0 &   \\
rc207 & OUR APPROACH & 841.5 & 857.2 & 1060.3 & 1000.0 & 309.2 &   \\
rc207 & BASELINE MILP & 602.0 & 738.3 & 988.1 & 1000.0 & 0.1 &   \\
rc208 & OUR APPROACH & 646.0 & 652.6 & 995.3 & 1000.0 & 324.3 &   \\
rc208 & BASELINE MILP & 536.8 & 609.2 & 854.5 & 1000.0 & 0.1 &   \\
\bottomrule
\end{tabular}
\end{table}
\begin{table}
\centering
\caption{Num Customers = 100 Problem Set = r100}
\label{tab_all:100_r100}
\begin{tabular}{cccccccccccccccccccccccccccccccccccccccccccccccccccccccc}
\toprule
file num & Approach & lp obj & mip dual bound & ILP obj & ilp time & total lp time & 10X+  \\
\midrule
r101 & OUR APPROACH & 1631.2 & 1637.7 & 1637.7 & 0.3 & 1.1 &   \\
r101 & BASELINE MILP & 1612.2 & 1637.7 & 1637.7 & 0.1 & 0.0 &   \\
r102 & OUR APPROACH & 1467.7 & 1467.7 & 1467.7 & 1.0 & 10.3 & YES \\
r102 & BASELINE MILP & 1027.2 & 1320.4 & 1470.5 & 1000.0 & 0.0 &    \\
r103 & OUR APPROACH & 1204.2 & 1208.7 & 1208.7 & 9.1 & 65.6 & YES \\
r103 & BASELINE MILP & 787.2 & 999.1 & 1244.2 & 1000.0 & 0.1 &    \\
r104 & OUR APPROACH & 944.9 & 946.5 & 1060.0 & 1000.0 & 338.2 &   \\
r104 & BASELINE MILP & 693.1 & 780.2 & 1012.8 & 1000.0 & 0.1 &   \\
r105 & OUR APPROACH & 1341.6 & 1355.3 & 1355.3 & 5.6 & 24.0 & YES \\
r105 & BASELINE MILP & 1098.2 & 1355.3 & 1355.3 & 440.9 & 0.0 &    \\
r106 & OUR APPROACH & 1215.4 & 1234.6 & 1234.6 & 70.2 & 95.6 & YES \\
r106 & BASELINE MILP & 833.3 & 1016.2 & 1266.9 & 1000.0 & 0.1 &    \\
r107 & OUR APPROACH & 1040.6 & 1049.8 & 1070.6 & 1000.0 & 201.7 &   \\
r107 & BASELINE MILP & 719.2 & 827.8 & 1126.8 & 1000.1 & 0.1 &   \\
r108 & OUR APPROACH & 906.9 & 908.6 & 1151.9 & 999.9 & 411.1 &   \\
r108 & BASELINE MILP & 679.8 & 736.3 & 997.1 & 1000.0 & 0.1 &   \\
r109 & OUR APPROACH & 1109.3 & 1136.7 & 1146.9 & 1000.0 & 109.2 &   \\
r109 & BASELINE MILP & 775.6 & 939.7 & 1194.9 & 1000.0 & 0.0 &   \\
r110 & OUR APPROACH & 1028.8 & 1034.9 & 1078.8 & 1000.2 & 188.2 &   \\
r110 & BASELINE MILP & 704.3 & 794.3 & 1115.9 & 1000.0 & 0.1 &   \\
r111 & OUR APPROACH & 1012.7 & 1026.3 & 1050.6 & 1000.0 & 193.2 &   \\
r111 & BASELINE MILP & 704.0 & 803.5 & 1126.2 & 1000.0 & 0.1 &   \\
r112 & OUR APPROACH & 910.0 & 911.1 & 1070.6 & 1000.0 & 438.7 &   \\
r112 & BASELINE MILP & 668.9 & 716.6 & 1069.9 & 1000.1 & 0.1 &   \\
\bottomrule
\end{tabular}
\end{table}
\begin{table}
\centering
\caption{Num Customers = 100 Problem Set = r200}
\label{tab_all:100_r200}
\begin{tabular}{cccccccccccccccccccccccccccccccccccccccccccccccccccccccc}
\toprule
file num & Approach & lp obj & mip dual bound & ILP obj & ilp time & total lp time & 10X+  \\
\midrule
r201 & OUR APPROACH & 1114.0 & 1143.2 & 1143.2 & 5.0 & 119.5 &   \\
r201 & BASELINE MILP & 979.2 & 1143.2 & 1143.2 & 21.2 & 0.0 &   \\
r202 & OUR APPROACH & 974.2 & 1006.9 & 1029.6 & 1000.0 & 362.9 &   \\
r202 & BASELINE MILP & 731.0 & 872.7 & 1049.5 & 1000.0 & 0.0 &   \\
r203 & OUR APPROACH & 794.0 & 809.0 & 883.2 & 1000.0 & 636.0 &   \\
r203 & BASELINE MILP & 633.0 & 716.7 & 906.3 & 1000.0 & 0.1 &   \\
r204 & OUR APPROACH & 668.0 & 675.8 & 770.1 & 1000.0 & 804.5 &   \\
r204 & BASELINE MILP & 591.1 & 647.9 & 789.8 & 1000.1 & 0.1 &   \\
r205 & OUR APPROACH & 891.0 & 912.4 & 954.7 & 1000.0 & 473.6 &   \\
r205 & BASELINE MILP & 738.1 & 879.8 & 967.4 & 1000.0 & 0.0 &   \\
r206 & OUR APPROACH & 803.7 & 825.6 & 898.0 & 1000.0 & 798.7 &   \\
r206 & BASELINE MILP & 651.9 & 758.2 & 914.6 & 1000.0 & 0.1 &   \\
r207 & OUR APPROACH & 719.4 & 728.3 & 900.0 & 1000.0 & 1687.1 &   \\
r207 & BASELINE MILP & 612.2 & 691.4 & 815.7 & 1000.1 & 0.1 &   \\
r208 & OUR APPROACH & 643.2 & 645.4 & 823.4 & 1000.0 & 396.1 &   \\
r208 & BASELINE MILP & 583.6 & 643.7 & 729.6 & 1000.0 & 0.1 &   \\
r209 & OUR APPROACH & 790.8 & 804.2 & 905.7 & 1000.0 & 339.0 &   \\
r209 & BASELINE MILP & 653.0 & 758.1 & 859.2 & 1000.0 & 0.1 &   \\
r210 & OUR APPROACH & 808.6 & 826.3 & 934.0 & 1000.0 & 407.7 &   \\
r210 & BASELINE MILP & 644.9 & 760.2 & 928.6 & 1000.0 & 0.1 &   \\
r211 & OUR APPROACH & 667.5 & 672.0 & 881.2 & 1000.0 & 292.8 &   \\
r211 & BASELINE MILP & 583.6 & 642.3 & 799.2 & 1000.1 & 0.1 &   \\
\bottomrule
\end{tabular}
\end{table}
\begin{table}
\centering
\caption{Num Customers = 100 Problem Set = c100}
\label{tab_all:100_c100}
\begin{tabular}{cccccccccccccccccccccccccccccccccccccccccccccccccccccccc}
\toprule
file num & Approach & lp obj & mip dual bound & ILP obj & ilp time & total lp time & 10X+  \\
\midrule
c101 & OUR APPROACH & 827.3 & 827.3 & 827.3 & 0.2 & 2.9 &   \\
c101 & BASELINE MILP & 819.7 & 827.3 & 827.3 & 0.1 & 0.1 &   \\
c102 & OUR APPROACH & 827.1 & 827.3 & 827.3 & 1.8 & 47.7 &   \\
c102 & BASELINE MILP & 754.3 & 827.3 & 827.3 & 24.1 & 0.0 &   \\
c103 & OUR APPROACH & 826.1 & 826.3 & 826.3 & 19.0 & 136.8 & YES \\
c103 & BASELINE MILP & 588.9 & 785.0 & 826.3 & 1000.0 & 0.1 &    \\
c104 & OUR APPROACH & 818.2 & 822.9 & 822.9 & 195.5 & 170.5 & YES \\
c104 & BASELINE MILP & 577.2 & 673.3 & 830.1 & 1000.0 & 0.1 &    \\
c105 & OUR APPROACH & 827.3 & 827.3 & 827.3 & 0.3 & 3.2 &   \\
c105 & BASELINE MILP & 818.2 & 827.3 & 827.3 & 0.1 & 0.0 &   \\
c106 & OUR APPROACH & 827.3 & 827.3 & 827.3 & 0.4 & 10.9 &   \\
c106 & BASELINE MILP & 699.7 & 827.3 & 827.3 & 0.7 & 0.0 &   \\
c107 & OUR APPROACH & 827.3 & 827.3 & 827.3 & 0.5 & 5.9 &   \\
c107 & BASELINE MILP & 818.1 & 827.3 & 827.3 & 0.4 & 0.0 &   \\
c108 & OUR APPROACH & 827.3 & 827.3 & 827.3 & 5.0 & 71.2 &   \\
c108 & BASELINE MILP & 569.5 & 827.3 & 827.3 & 8.4 & 0.0 &   \\
c109 & OUR APPROACH & 825.1 & 827.3 & 827.3 & 54.8 & 145.2 & YES \\
c109 & BASELINE MILP & 566.0 & 792.1 & 853.4 & 1000.0 & 0.1 &    \\
\bottomrule
\end{tabular}
\end{table}
\begin{table}
\centering
\caption{Num Customers = 100 Problem Set = c200}
\label{tab_all:100_c200}
\begin{tabular}{cccccccccccccccccccccccccccccccccccccccccccccccccccccccc}
\toprule
file num & Approach & lp obj & mip dual bound & ILP obj & ilp time & total lp time & 10X+  \\
\midrule
c201 & OUR APPROACH & 589.1 & 589.1 & 589.1 & 0.1 & 5.5 &   \\
c201 & BASELINE MILP & 589.1 & 589.1 & 589.1 & 0.1 & 0.0 &   \\
c202 & OUR APPROACH & 589.1 & 589.1 & 589.1 & 0.3 & 65.6 &   \\
c202 & BASELINE MILP & 556.7 & 589.1 & 589.1 & 1.6 & 0.1 &   \\
c203 & OUR APPROACH & 583.7 & 588.7 & 588.7 & 13.9 & 127.4 &   \\
c203 & BASELINE MILP & 534.9 & 588.7 & 588.7 & 6.0 & 0.1 &   \\
c204 & OUR APPROACH & 582.4 & 588.1 & 588.1 & 61.6 & 183.7 &   \\
c204 & BASELINE MILP & 517.5 & 588.1 & 588.1 & 77.8 & 0.1 &   \\
c205 & OUR APPROACH & 586.0 & 586.4 & 586.4 & 0.3 & 29.5 &   \\
c205 & BASELINE MILP & 539.7 & 586.4 & 586.4 & 0.4 & 0.0 &   \\
c206 & OUR APPROACH & 586.0 & 586.0 & 586.0 & 0.2 & 33.1 &   \\
c206 & BASELINE MILP & 536.3 & 586.0 & 586.0 & 2.1 & 0.0 &   \\
c207 & OUR APPROACH & 582.7 & 585.8 & 585.8 & 5.2 & 98.5 &   \\
c207 & BASELINE MILP & 535.7 & 585.8 & 585.8 & 2.0 & 0.0 &   \\
c208 & OUR APPROACH & 583.4 & 585.8 & 585.8 & 1.0 & 72.0 &   \\
c208 & BASELINE MILP & 531.9 & 585.8 & 585.8 & 3.3 & 0.0 &   \\
\bottomrule
\end{tabular}
\end{table}
\end{document}